\documentclass[12pt]{article}
\usepackage{amssymb,amsmath}
\usepackage{cite}
\usepackage{arydshln}

\newtheorem{lemma}{Lemma}
\newtheorem{prop}{Proposition}
\newtheorem{thm}{Theorem}
\newtheorem{cor}{Corollary}

\def\id{\mbox{id}}

\begin{document}

\begin{center}
{\large\bf Invariants and reduced Wigner coefficients for quasi-triangular Hopf superalgebras}\\
~~\\

{\large Mark D. Gould, Phillip S. Isaac and Jason L. Werry}\\
~~\\

School of Mathematics and Physics, The University of Queensland, St Lucia QLD 4072, Australia.
\end{center}

\begin{abstract}
We develop explicit formulae for the eigenvalues of various invariants 
for highest weight irreducible representations of the quantum supergroup $U_q[gl(m|n)]$. The
techniques employed make use of modified characteristic identity methods and allow for the
evaluation of generator matrix elements and reduced Wigner coefficients.
\end{abstract}

\section{Introduction}

In a 2019 article by the first two authors \cite{GI2019}, explicit formulae for an assortment of
invariants and generator matrix elements were constructed using charateristic identity techniques
specific to the case of $U_q[gl(n,\mathbb{C})].$ In a sense, this paper 
complimented certain earlier works on quantum groups, for example \cite{GLB1992,Gould1992}. In that article
\cite{GI2019}, it was argued that the standard
characteristic matrix used in the construction for the cases of Lie algebras and superalgebras, could
not be utilised in the quantum group case due to block partitioning issues. The main motivation for
\cite{GI2019} was to resolve the issue in the (non-graded) quantum group case, and to this end the
$R$-matrix was gainfully employed. For further background and motivation to the broader
problem we direct the reader to that work \cite{GI2019}.

This research sits within the realm of {\em constructive representation theory}, inspired by the work of
Baird and Biedenharn \cite{BB1963}, the aim of which is to
establish explicit and accessible formulae for significant representation-theoretic quantities such as matrix elements
and Wigner coefficients. Explicit expressions for such quantities can be useful, for example, in developing a
systematic approach to the construction of solutions of the Yang-Baxter equation (with spectral
parameter) \cite{R1988} which are of importance to quantum integrable systems and associated link invariants. 
Within this area of research, powerful techniques utilising characteristic identities were introduced
by Green and Bracken \cite{Green1971,BraGre1971}, and largely developed in the early work of Gould
\cite{Gould1978,Gould1980,Gould1981,Gould1981b,Gould1986,Gould1986b,Gould1987}. In more recent times, similar
formulae have been further developed by the current authors for the general linear and orthosymplectic Lie
superalgebras \cite{GIW1,GIW2,GIW3,GIW4,GIW5,GI2015}.

As promised in that 2019 article, in this current work we develop analogous formulae for the quantum
supergroup $U_q[gl(m|n)]$. With care, we do indeed see that the techniques making use of the $R$-matrix in
solving the block partitioning problem of the characteristic matrix extend to this
$\mathbb{Z}_2$-graded case.


After introducing some preliminary notation and discussing key objects such as tensor operators in
Section \ref{Section2}, we develop the modified characteristic matrix in Section \ref{Section3} for
the case of $U_q[gl(m|n)]$. This then allows us to develop formulae for the eigenvalues of certain
invariants, which is done in the concluding Section \ref{Section4}, with an emphasis on the reduced Wigner
coefficients.


\section{Quasi-triangular Hopf superalgebras} \label{Section2}

Let $H$ be a (quasi-triangular) Hopf superalgebra with product $m$, co-product $\Delta$, co-unit $\varepsilon$
and antipode $S$. Then $H$ is also a (quasi-triangular) Hopf superalgebra with opposite co-product
$\Delta^T = T\circ\Delta$, with $T$ the linear twist operator, and antipode $S^{-1}$. Moreover, if $R$ is the
$R$-matrix then $R^T= T(R)$ is the $R$-matrix under this opposite structure.

Throughout we use the notation of Sweedler \cite{S69,M95} for the co-product, so that
$$
\Delta(a) = a_{(1)}\otimes a_{(2)},
$$
where there is an implicit summation that depends on the element $a\in H$.
We describe the $\mathbb{Z}_2$-grading using the following notation. For $H=H_0\oplus H_1$, we call
elements of $H_0\cup H_1$ homogeneous, and write $[a] = \alpha$, for homogeneous $a\in H_\alpha$
($\alpha=0,1$). In what follows, we assume an element $a$ is homogeneous when we write $[a]$. 
Recall the linear twist operator on $H\otimes H$, defined on homogeneous elements as
$$
T(a\otimes b) = (-1)^{[a][b]}b\otimes a,
$$
so
$$
\Delta^T(a) = (-1)^{[a_{(1)}][a_{(2)}]}a_{(2)}\otimes a_{(1)}.
$$
Thus
$$
S^{-1}(a_{(2)})a_{(1)}(-1)^{[a_{(1)}][a_{(2)}]} = (-1)^{[a_{(1)}][a_{(2)}]}a_{(2)}S^{-1}(a_{(1)}) =
\varepsilon(a).
$$
Also recall the anti-homomorphism property of $S$ so
$$
S(ab) = (-1)^{[a][b]}S(b)S(a).
$$
Finally, if we write
$$
R = a_i\otimes b_i \mbox{ (sum on $i$)} 
$$
then
$$
R^T = T(R) = (-1)^{[a_i]}b_i\otimes a_i.
$$
\noindent
\underline{{\em Remark}:} It is worth pointing out that $R$ is an even element of $H\otimes H$, and
so the grading of tensor product components of $R=a_i\otimes b_i$ satisfies
$[a_i]=[b_i]$ for each $i$. Therefore, $R^T=(-1)^{[a_i][b_i]}b_i\otimes a_i=(-1)^{[a_i]}b_i\otimes a_i$, for example.

By definition, $R$ satisfies
\begin{align}
(\Delta\otimes\id)R &= R_{13}R_{23}, \label{p1stara} \\
(\id\otimes\Delta)R &= R_{13}R_{12}, \label{p1starb}
\end{align}
and the intertwining property
$$
R\Delta(a) = \Delta^T(a)R,
$$
so
\begin{equation}
(-1)^{[b_i][a_{(1)}]}a_ia_{(1)}\otimes b_ia_{(2)} = (-1)^{[a_{(1)}]([a_{(2)}] +
[a_i])}a_{(2)}a_i\otimes a_{(1)}b_i.
\label{p1starstar}
\end{equation}

Applying $(\varepsilon\otimes\id\otimes\id)$ and $(\id\otimes\varepsilon\otimes\id)$ to
(\ref{p1stara}) and (\ref{p1starb}) gives
\begin{equation}
(\varepsilon\otimes\id)R = (\id\otimes\varepsilon)R = I. \label{rcounit}
\end{equation}
Then multiply $(\ref{p1stara})$ by $(S\otimes\id\otimes\id)$ or $(\id\otimes S\otimes \id)$ and
apply $m\otimes \id$ to give
$$
I\otimes I = (S\otimes \id)R\cdot R \ \ \Rightarrow \ \ R^{-1} = (S\otimes\id)R = (\id\otimes
S^{-1})R.
$$
\underline{{\em Remark}:} It is worth recalling the antipode property:
$$
\varepsilon(a)I = a_{(1)}S(a_{(2)}) = S(a_{(1)})a_{(2)}.
$$ 
Note, however, for $S^{-1}$ we need to include the grading factor $(-1)^{[a_{(1)}][a_{(2)}]}$ which
arises in the definition of $\Delta^T$.

Following the non-graded case, we consider the $u$-operator \cite{D90,R1990}. With $R=a_i\otimes
b_i$ (sum on $i$), we define
\begin{equation}
u = m\circ (S\otimes \id)\circ T(R).
\label{uoperator}
\end{equation}
The properties of $u$ that we will make use of in this article are as follows. We note that the
proofs of these results are presented in Appendix \ref{AppendixA}.
For $u$ defined as in equation (\ref{uoperator}) above, we have

\ \\

(i) $S^2(a)u = ua,$ $\forall a\in H,$

\ \\

(ii) $u$ is invertible with $u^{-1} = m\circ(S^{-1}\otimes\id)\circ T(R^{-1}),$ 

\ \\

(iii) $S^2(a) = uau^{-1},\ \ \forall a\in H,$ 

\ \\

(iv) $\Delta(u) = (u\otimes u)(R^TR)^{-1}.$

\ \\

\noindent
We further observe that 
$
\tilde{R} = (R^T)^{-1}
$
is also an $R$-matrix. Denoting the corresponding
$u$-operator by $\tilde{u}$, we then have 

\ \\

(v) $\tilde{u} = S(u^{-1}).$

\ \\

\noindent
Finally,

\ \\

(vi) $u\otimes u$ commutes with $R^TR$,

\ \\

(vii) $\Delta(S(u)) = (S(u)\otimes S(u))(R^TR)^{-1}$.


\subsection{Tensor operators}

Let $V$ be a finite-dimensional irreducible $H$-module and $\pi$ the representation afforded by $V$.
Given a basis $\{e_\alpha\}_{\alpha=1}^n$ for $V$ we define a collection of operators
$$
T\equiv \{ T_{\alpha}\}_{\alpha=1}^n,
$$
transforming according to the rule
\begin{equation}
(-1)^{[a][\tau]}aT_\alpha = T_\beta \pi(a_{(1)})_{\beta\alpha}a_{(2)}(-1)^{[a_{(2)}][\alpha]},
\label{p10star}
\end{equation}
to be an irreducible tensor operator of type (or rank) $\pi$. 
Here $[\tau]$ denotes the parity of the tensor operator so $[T_{\alpha}] = [\tau] + [\alpha],$ where
we use $[\alpha]$ to denote the grading of the basis element $e_\alpha$, i.e.
$[\alpha] \equiv [e_\alpha].$

\noindent
\underline{{\em Remark}:} The usual definition assumes $[\tau]=0$ corresponding to an even tensor
operator. Here we allow also for the possibility that $[\tau]=1$, corresponding to an odd tensor.

\begin{lemma}
Let $W$ be an $H$-module and set 
$$
X = \mbox{span}\{T_\alpha w|w\in W, \alpha=1,\ldots,n\}
$$
(also an $H$-module). Then $T:V\otimes W\rightarrow X,$ $T(e_\alpha\otimes w) = T_\alpha w$ is an
intertwining operator, i.e.
$$
a\circ T = (-1)^{[a][\tau]}T\Delta(a), \ \ a\in H.
$$
\end{lemma}
Proof:
\begin{align*}
T\Delta(a)(e_\alpha\otimes w) &= T(a_{(1)}e_\alpha\otimes a_{(2)}w)(-1)^{[a_{(2)}][\alpha]}\\
&= \pi(a_{(1)})_{\beta \alpha} T(e_\beta\otimes a_{(2)}w)(-1)^{[a_{(2)}][\alpha]}\\
&=  \pi(a_{(1)})_{\beta \alpha} T_\beta a_{(2)}w (-1)^{[a_{(2)}][\alpha]}\\
& \stackrel{(\ref{p10star})}{=} (-1)^{[a][\tau]}a\cdot T_\alpha w\\
&= (-1)^{[a][\tau]}a \cdot T(e_\alpha\otimes w).
\end{align*}

\begin{flushright}$\Box$\end{flushright}

\noindent
\underline{{\em Remarks}:} It is suggestive that we define the ``adjoint'' action in the
usual way by 
$$
\mbox{Ad}a\circ T_\alpha = a_{(1)} T_\alpha S(a_{(2)}) (-1)^{[a_{(2)}]([\tau]+[\alpha])}.
$$
Then we have 
\begin{align*}
\mbox{Ad}a\circ T_\alpha 
&= 
(a_{(1)}T_\alpha)\cdot\left\{ S(a_{(2)}) (-1)^{[a_{(2)}]([\tau]+[\alpha])} \right\} \\
&\stackrel{(\ref{p10star})}{=}
\left\{ (-1)^{([a_{(1)}]+[a_{(2)}])[\tau]} \pi(a_{(1)})_{\beta \alpha} T_\beta a_{(2)}
(-1)^{[a_{(2)}][\alpha]} \right\} S(a_{(3)}) (-1)^{[a_{(3)}]([\tau]+[\alpha])} \\
&=
(-1)^{[a][\tau]} \pi(a_{(1)})_{\beta\alpha} T_\beta a_{(2)}S(a_{(3)})
(-1)^{[\alpha]([a_{(2)}]+[a_{(3)}])} \\
&= 
(-1)^{[a][\tau]} \pi(a)_{\beta\alpha}T_\beta
\end{align*}
which we may also take as the transformation law for tensor operators. We prefer, however,
the definition (\ref{p10star}).

\noindent
\underline{{\em Example}:} With the notation as above we may consider
$$
T = (\pi\otimes \id)R^TR = \sum_{\alpha,\beta} e_{\alpha \beta}\otimes T_{\alpha\beta}
$$
(this defines $T_{\alpha\beta}$). We then have
\begin{align*}
T(e_\alpha\otimes w) &= (e_{\beta\delta}\otimes T_{\beta\delta})(e_\alpha\otimes w)\\
&=
(-1)^{([\beta]+[\delta])[\alpha]} e_{\beta\delta}e_\alpha\otimes T_{\beta\delta}w\\
&=
(-1)^{[\alpha]([\alpha]+[\beta])} e_\beta\otimes T_{\beta\alpha}w.
\end{align*}
As we shall see, this determines a tensor operator.


\subsection{Basic module constructions}\label{Section5} 

Given $H$-modules $V,W$, recall that $\ell(V,W)$, the space of linear maps from $V$ to $W$, is an
$H$-module under the action defined by
$$
a\circ\varphi(v) = (-1)^{[\varphi][a_{(2)}]}a_{(1)}\varphi(S(a_{(2)})v),\ \ \forall v\in V,\
\varphi\in \ell(V,W).
$$
Then $\varphi$ is an $H$-module homomorphism if and only if it is an even invariant under this
action, i.e.
$$
a_{(1)}\varphi(S(a_{2})v) = \varepsilon(a)\varphi(v).
$$
In the case $W=\mathbb{C}$ is the trivial 1-dimensional module, this action reduces to
\begin{align*}
a\circ\varphi(v) &= \varepsilon(a_{(1)})\varphi(S(a_{(2)})v)(-1)^{[\varphi][a_{(2)}]}\\
&= (-1)^{[\varphi][a]}\varphi(S(a)v),\ \forall \varphi\in V^*,v\in V.
\end{align*}
\underline{{\em Note}}: Let $\{ e_\alpha\}$ be a basis for the finite-dimensional $H$-module $V$
with dual basis $\{e_\alpha^*\}$ for $V^*$ defined by
$$
\langle e_\alpha^*,e_\beta \rangle \equiv e_\alpha^*(e_\beta) = \delta_{\alpha\beta}.
$$
Then for $a\in H$ we may write
$$
ae_{\alpha} = \langle e_\beta^*,ae_\alpha \rangle e_\beta \mbox{ (sum on $\beta$)}
$$
and similarly
$$
ae_{\alpha}^* = \langle ae_\alpha^*,e_\beta\rangle e_\beta^* \mbox{ (sum on $\beta$)}
$$
where we now define, as above,
$$
\langle ae_\alpha^*,e_\beta\rangle = (-1)^{[\alpha][a]}\langle e_\alpha^*,S(a)e_\beta\rangle.
$$
Then if $\pi$, $\pi^*$ are the representations afforded by $V$, $V^*$ respectively, we have the
matrix elements
\begin{align*}
\pi(a)_{\beta\alpha} &= \langle e_\beta^*,ae_\alpha \rangle \\
\pi^*(a)_{\beta\alpha} &= \langle ae_\alpha^*,e_\beta \rangle \\
&= \langle e_\alpha^*,S(a)e_\beta\rangle (-1)^{[\alpha][a]}\\
&= (-1)^{[\alpha][a]} \pi(S(a))_{\alpha\beta}\\
\Rightarrow \ \ \pi^*(a) &= \pi(S(a))^T
\end{align*}
where $T$ is the ``super transpose'' defined by
$$
\left( A^T \right)_{\alpha\beta} = (-1)^{[\beta][A]}A_{\beta\alpha}.
$$

\begin{prop}
Let $V,W$ be finite dimensional $H$-modules. Then we have the $H$-module isomorphism
$$
W\otimes V^* \cong \ell(V,W).
$$
\end{prop}
Proof:
Following the classical argument, define a linear map $\varphi:W\otimes V^*\longrightarrow
\ell(V,W)$ by
$$
\varphi(w\otimes v^*)(u) = v^*(u)w,\ \forall u\in V,\ w\in W,\ v^*\in V^*.
$$
Then $\varphi$ is well-defined and one-to-one. It is also onto since given $f\in \ell(V,W)$ we have
$$
f = \varphi(f(e_\alpha)\otimes e^*_\alpha) \ \ \mbox{(sum on $\alpha$)},
$$
where, as above, $\{e_\alpha\}$ is a basis for $V$ with dual basis $\{ e_\alpha^*\}$ for $V^*$. It
remains to check that $\varphi$ is a $H$-module homomorphism.

To this end we have
\begin{align*}
\varphi(\Delta(a)(w\otimes v^*))(u) &= (-1)^{[a_{(2)}][w]}\varphi(a_{(1)}w\otimes a_{(2)}v^*)(u)\\
&= (-1)^{[a_{(2)}][w]}\langle a_{(2)}v^*,u \rangle a_{(1)}w\\
&= (-1)^{[a_{(2)}]([w]+[v^*])}\langle v^*,S(a_{(2)})u \rangle a_{(1)}w\\
&= a_{(1)}\left\{ \varphi(w\otimes v^*)(S(a_{(2)})u) \right\} (-1)^{[a_{(2)}]([w]+[v^*])} \\
&= (a\circ \varphi(w\otimes v^*))(u),\ \ \forall u\\
\Rightarrow \ \ \varphi(\Delta(a)(w\otimes v^*)) &= a\circ\varphi(w\otimes v^*).
\end{align*}
\begin{flushright}$\Box$\end{flushright}

In the case $W=V$ is finite dimensional and irreducible we arrive at the following.

\begin{lemma} \label{lem14}
The identity module occurs exactly once in $V\otimes V^*$ and is spanned by the vector (notation as
above)
$$
v_0=e_\alpha\otimes e_\alpha^*\ \ \mbox{(sum on $\alpha$)}.
$$
\end{lemma}
Proof: Since $V\otimes V^*\cong \ell(V,V)$, Schur's Lemma implies the identity module occurs exactly
once. By direct calculation we have for $a\in H$
\begin{align*}
a\circ v_0 &= \Delta(a)(e_\alpha\otimes e_\alpha^*)\\
&= (-1)^{[a_{(2)}][\alpha]} a_{(1)}e_\alpha\otimes a_{(2)}e_\alpha^*\\
&= \langle a_{(2)}e_\alpha^*,e_\beta \rangle (-1)^{[a_{(2)}][\alpha]}a_{(1)}e_\alpha\otimes
e_\beta^*\\
&= \langle e_\alpha^*,S(a_{(2)})e_\beta \rangle a_{(1)}e_\alpha\otimes e_\beta^*\\
&= a_{(1)}S(a_{(2)})e_\beta\otimes e_\beta^*\\
&= \varepsilon(a)v_0. 
\end{align*}
\begin{flushright}$\Box$\end{flushright}

\underline{{\em Remark}:} The identity module also occurs in $V^*\otimes V$ and is spanned by the
vector
$$
\Omega = (-1)^{[\alpha]}e_\alpha^*\otimes u^{-1}e_\alpha.
$$
Indeed we have
\begin{align*}
\Delta(a)\Omega &= (-1)^{[\alpha]+[\alpha][a_{(2)}]}a_{(1)}e_\alpha^*\otimes a_{(2)}u^{-1}e_\alpha\\
&= (-1)^{[\alpha]+[\alpha][a_{(2)}]} \langle a_{(1)}e_\alpha^*,e_\beta \rangle e_\beta^*\otimes
a_{(2)}u^{-1}e_\alpha\\
&= (-1)^{[\alpha]+[\alpha][a]} \langle e_\alpha^*,S(a_{(1)})e_\beta \rangle e_\beta^*\otimes
a_{(2)}u^{-1}e_\alpha\\
&= (-1)^{[\alpha]+[\alpha][a]} e_\beta^*\otimes a_{(2)}u^{-1}\langle e_\alpha^*,S(a_{(1)})e_\beta
\rangle e_\alpha.
\end{align*}
Now observe that
\begin{align*}
[\alpha]+[\alpha][a] &= [a_{(1)}]+ [\beta] + [a_{(1)}][a]+[\beta][a]\\
&= [a_{(1)}]+ [\beta] + [a_{(1)}]+  [a_{(1)}] [a_{(2)}]+[\beta][a]\\
&= [\beta] +[\beta][a]+  [a_{(1)}] [a_{(2)}]
\end{align*}
\begin{align*}
\Rightarrow \ \ \Delta(a)\Omega &= (-1)^{[\beta] +[\beta][a]+  [a_{(1)}] [a_{(2)}]} e_\beta^*\otimes
a_{(2)}u^{-1}S(a_{(1)})e_\beta\\
&= (-1)^{[\beta] +[\beta][a]} e_\beta^*\otimes a_{(2)}S^{-1}(a_{(1)})u^{-1}e_\beta(-1)^{[a_{(1)}]
[a_{(2)}]} \\
&= (-1)^{[\beta] +[\beta][a]} e_\beta^*\otimes \varepsilon(a) u^{-1}e_\beta\\
&= \varepsilon(a) (-1)^{[\beta]} e_\beta^*\otimes u^{-1}e_\beta\\
&= \varepsilon(a)\Omega.
\end{align*}
The above invariants in $V\otimes V^*$ and $V^*\otimes V$ are of importance below.


\subsection{Tensor operators from $R^TR$}\label{Section6} 

As above, let $\pi$ be the representation afforded by a finite-dimensional irreducible $H$-module
$V$. Then $T = (\pi\otimes\id)(R^TR)$ acting on an $H$-module $W$ determines a matrix with entries
$T_{\alpha\beta}$ defined by 
\begin{align*}
T(e_\beta\otimes w) &= (\pi\otimes\id)(R^TR)(e_\beta\otimes w)\\
&= e_\alpha\otimes T_{\alpha\beta}w \mbox{ (sum on $\alpha$), } w\in W. 
\end{align*}
The $T_{\alpha\beta}$ transform as a tensor operator since
\begin{prop} \label{prop16}
$\varphi: V^*\otimes V\otimes W\longrightarrow W$, defined by
\begin{align*}
\varphi(e_\alpha^*\otimes e_\beta\otimes w) &= (e_\alpha^*\otimes I)[(\pi\otimes
\id)(R^TR)](e_\beta\otimes w)\\
&= (e_\alpha^*\otimes I)T(e_\beta\otimes w) = T_{\alpha\beta}w
\end{align*}
is an intertwining operator.
\end{prop}
Proof: 
\begin{align*}
\varphi((\id\otimes\Delta)\Delta(a)\cdot(e_\alpha^*\otimes e_\beta\otimes w) )
&= \varphi(a_{(1)}e_\alpha^*\otimes a_{(2)}e_\beta\otimes a_{(3)}w) (-1)^{[a_{(2)}][\alpha] +
[a_{(3)}]([\alpha]+[\beta])}\\
&= (a_{(1)}e_\alpha^*\otimes I)\circ(R^TR)(a_{(2)}e_\beta\otimes a_{(3)}w)(-1)^{[a_{(2)}][\alpha] +
[a_{(3)}]([\alpha]+[\beta])}\\
&= (a_{(1)}e_\alpha^*\otimes I)\circ(R^TR)\Delta(a_{(2)})(e_\beta\otimes
w)(-1)^{[a_{(2)}][\alpha]}\\
&= (e_\alpha^*\otimes I)(S(a_{(1)})\otimes I)\Delta(a_{(2)})R^TR(e_\beta\otimes w)(-1)^{[a_{(2)}][\alpha]}\\
&= (e_\alpha^*\otimes I)(I\otimes a)R^TR (e_\beta\otimes w)(-1)^{[a][\alpha]}\\
&= a\cdot \varphi(e_\alpha^*\otimes e_\beta\otimes w).
\end{align*}

\begin{flushright}$\Box$\end{flushright}

The above transformation law gives
\begin{align*}
a T_{\alpha\beta}w &= \varphi((\id\otimes \Delta)\Delta(a)(e_\alpha^*\otimes
e_\beta\otimes w))\\
&=
(-1)^{[a_{(2)}][\alpha]+[a_{(3)}]([\alpha]+[\beta])}
\varphi(a_{(1)}e_\alpha^*\otimes a_{(2)}e_\beta \otimes a_{(3)}w) \\
&=
(-1)^{[a_{(2)}][\alpha]+[a_{(3)}]([\alpha]+[\beta])}
\pi^*(a_{(1)})_{\gamma\alpha}\pi(a_{(2)})_{\delta\beta}T_{\gamma\delta}a_{(3)}w \\
&=
(\pi^*\otimes\pi)(\Delta(a_{(1)}))_{\gamma\delta,\alpha\beta} T_{\gamma\delta}a_{(2)}w
(-1)^{[a_{(2)}]([\alpha]+[\beta])} 
\end{align*}
which is the required transformation law. We write this as
$$
aT_{\alpha\beta} = (-1)^{[a_{(2)}][\alpha]+[a_{(3)}]([\alpha]+[\beta])}
\pi^*(a_{(1)})_{\gamma\alpha}\pi(a_{(2)})_{\delta\beta} T_{\gamma\delta} a_{(3)}.
$$



\section{Quantum supergroups and characteristic matrix for $U_q[gl(m|n)]$}\label{Section3} 

Let $L$ be a simple Lie superalgebra and $U_q(L)$ the corresponding quantum supergroup with simple
generators $e_i$, $f_i$, $h_i$ ($1\leq i\leq \ell=$ rank$L$), co-unit
$$
\varepsilon(x) = 0,\ x=e_i,f_i,h_i,\ \ \varepsilon(I)=1,
$$
co-product
\begin{align*}
\Delta(x) &= q^{h_i/2}\otimes x + x\otimes q^{-h_i/2},\ \ x=e_i,f_i,\\
\Delta(h) &= h\otimes I + I\otimes h,\ \ h\in H \mbox{ (Cartan subalgebra),}
\end{align*}
and antipode
$$
S(I)=I,\ \ S(x) = -q^{-h_i/2}xq^{h_i/2} = -q^{-h_\rho}xq^{h_\rho},\ \ x=e_i,f_i,h
$$
where $\rho$ is the graded half sum of the positive roots
$$
\rho = \rho_0-\rho_1,\ \ \rho_0 = \frac12\sum_{\alpha\in\Phi_0^+}\alpha,\ \ 
\rho_1 = \frac12\sum_{\beta\in\Phi_1^+}\beta,
$$
where we have denoted the sets of even and odd positive roots as $\Phi_0^+$ and $\Phi_1^+$
respectively.
\underline{{\em Note}}: For $L=gl(m|n)$ we have in usual notation
\begin{align}
h_a &= (-1)^{[a]}(E_{aa} - (-1)^{[\alpha_a]}E_{a+1\, a+1}),\ \ 1\leq a<m+n, \nonumber\\
\rho &= \frac12\sum_{i=1}^m(m-n-2i+1)\varepsilon_i +
\frac12\sum_{\mu=1}^n(m+n-2\mu+1)\delta_\mu, \label{p18star}
\end{align}
where $h_a$ is defined by
$$
\lambda(h_a) = (\lambda,\alpha_a),\ \ \lambda\in H^*,
$$
with the invariant form $(\, ,\, )$ on $H^*$ defined by
$$
(\varepsilon_i,\varepsilon_j)=\delta_{ij},\ \ (\varepsilon_i,\delta_\mu)=0,\ \
(\delta_\mu,\delta_\nu)=-\delta_{\mu\nu}.
$$

The universal $R$-matrix for $U_q(L)$ is expressible
\begin{equation}
R = q^{h_i\otimes h^i}\left\{ I\otimes I + \sum_s e_s\otimes e^s \right\} \in U_q^{(-)}(L)\otimes
U_q^{(+)}(L)
\label{p19star}
\end{equation}
with $\{h_i|i=1,\ldots,\ell\}$ a basis for $H$ with dual basis $\{h^i|i=1,\ldots,\ell\}$
under the invariant bilinear form $(,)$ on $L$ so that
$$
\lambda(h_i)\mu(h^i) = (\lambda,\mu),\ \ \forall \lambda,\mu\in H^*.
$$
Also, $U_q^{(\pm)}(L)$ are the Hopf subalgebras generated by $H$ (or more precisely
$q^{h_i}$) together with the raising simple generators $e_i$ (respectively the lowering
generators $f_i$).

It follows in this case that 
\begin{align*}
S^2(a) &= q^{-2h_\rho} a q^{2h_\rho} , \ \ \forall a\in L\\
&= uau^{-1}\\
\Rightarrow \ \ v &= q^{2h_\rho} u
\end{align*}
is an even central element with inverse
$$
v^{-1} = u^{-1}q^{-2h_\rho}.
$$
Moreover, since the $q^{\pm h_\rho}$ are group-like, we also have 
$$
\Delta(v) = (v\otimes v)(R^TR)^{-1}
$$
$$
\Rightarrow \ \ R^TR = (v\otimes v)\Delta(v^{-1}).
$$
Similarly utilising the $R$-matrix $\tilde{R}=(R^T)^{-1}$,
we have, in terms of the corresponding $u$-operator $\tilde{u} = S(u^{-1}),$ the central
element
\begin{align*}
\tilde{v} &= q^{2h_\rho} \tilde{u}\\
&= q^{2h_\rho}S(u^{-1})\\
&= S(u^{-1}q^{-2h_\rho}) = S(v^{-1}).
\end{align*}
In this case
\begin{align*}
\Delta(\tilde{v}) &= (\tilde{v}\otimes\tilde{v})(\tilde{R}^T\tilde{R})^{-1}\\
&= (\tilde{v}\otimes\tilde{v}) (R^TR) = S(v^{-1})\otimes S(v^{-1}) (R^TR).
\end{align*}

Regarding eigenvalues, first observe that with $R$ as in (\ref{p19star}) that
$\tilde{R}=(R^T)^{-1}$ has the form
$$
\tilde{R} = q^{-h_i\otimes h^i}\left\{ I\otimes I + \sum_s e^s\otimes e_s \right\} \in
U_q^{(+)}(L)\otimes U_q^{(-)}(L).
$$
It follows that if $V(\Lambda)$ is a finite dimensional
irreducible $U_q(L)$-module with highest weight $\Lambda$ and maximal weight vector
$e_+^\Lambda$ that
\begin{align*}
\tilde{u}e_+^\Lambda &= q^{h_ih^i}e_+^\Lambda = q^{(\Lambda,\Lambda)}e_+^\Lambda\\
\Rightarrow \ \ \tilde{v}e_+^\Lambda &= q^{2h_\rho} \tilde{u} e_+^\Lambda =
q^{(\Lambda,\Lambda+2\rho)}e_+^\Lambda. 
\end{align*}
Therefore the eigenvalue of $\tilde{v}$ in $V(\Lambda)$ is given by
$$
\chi_\Lambda(\tilde{v}) = q^{(\Lambda,\Lambda+2\rho)},
$$
and so $v=S^{-1}(\tilde{v}^{-1})$ takes the eigenvalue 
$$
\chi_\Lambda(v) = q^{-(\Lambda,\Lambda+2\rho)}.
$$
\underline{{\em Remark}:}  Let $U_q(L_{\bar{0}})$ be the quantum group corresponding to
the even Lie subalgebra of $L=L_{\bar{0}}\oplus L_{\bar{1}}$ ($\mathbb{Z}_2$
decomposition). Then an alternative way of evaluating the eigenvalue of $v$ directly is to
consider the Kac module $K(\Lambda)$ which has minimal weight $\Lambda_- - 2\rho_1$ where
$\Lambda_-$ is the minimal weight of the irreducible $U_q(L_{\bar{0}})$-module
$V_{\bar{0}}(\Lambda)$. Then $\Lambda_-=\tau(\Lambda)$ where $\tau$ is the unique element
of the Weyl group of $L_{\bar{0}}$ sending even positive roots to negative ones. Then in
view of the form (\ref{p19star}) of the $R$-matrix we have on the minimal weight vector
$e_-^\Lambda\in K(\Lambda),$
\begin{align*}
ue_-^\Lambda &= q^{-(\Lambda_--2\rho_1,\Lambda_--2\rho_1)}e_-^\Lambda\\
\Rightarrow \ \ ve_-^{\Lambda} &= q^{2h_\rho}ue_-^\Lambda = 
q^{-(\Lambda_--2\rho_1,\Lambda_--2\rho_1)+2(\rho_0-\rho_1,\Lambda_--2\rho_1)}e_-^\Lambda.
\end{align*}
Using $\tau(\rho_0) = -\rho_0$ and $\tau(\rho_1) = \rho_1$, we obtain, using Weyl group
invariance of the form, for the above $q$-exponent,
\begin{align*}
-(\Lambda_--2\rho_1,\Lambda_--2\rho_1)+2(\rho_0-\rho_1,\Lambda_--2\rho_1)
&= -(\Lambda_--2\rho_1,\Lambda_--2\rho_1+2\rho_1-2\rho_0)\\
&= -(\Lambda_--2\rho_1,\Lambda_--2\rho_0)\\
&= -(\Lambda_-, \Lambda_--2\rho_1-2\rho_0)-4(\rho_1,\rho_0)\\
&= -(\Lambda_-,\Lambda_--2\rho_1-2\rho_0) \mbox{ (since $(\rho_1,\rho_0)=0$)}\\
&= -(\tau(\Lambda_-),\tau(\Lambda_-)-2\tau(\rho_1) - 2\tau(\rho_0) )\\
&= -(\Lambda,\Lambda-2\rho_1+2\rho_0)\\
&= -(\Lambda,\Lambda+2\rho)
\end{align*}
$$
\Rightarrow \ \ ve_-^\Lambda = q^{-(\Lambda,\Lambda+2\rho)}e_-^\Lambda,
$$
as expected.


\subsection{Characteristic identities}\label{Section8} 

Let $V(\Lambda)$ be a finite-dimensional irreducible $U_q(L)$-module with highest weight $\Lambda$
and set
$$
A = (q-q^{-1})^{-1}(\pi_\Lambda\otimes\id)(I\otimes I-R^TR),
$$
where $\pi_\Lambda$ is the irreducible representation afforded by $V(\Lambda)$. In the classical
limit $q\rightarrow 1$ this gives the characteristic matrix considered previously. 
Since $R^TR$ transforms as a tensor operator, the entries
$A_{\alpha\beta}$ of the matrix $A$ also transform as a tensor operator of type
$\pi_\Lambda^*\otimes\pi_\Lambda$.

Acting on a $U_q(L)$-module $W$, however, the matrix $A$ is an operator on $V(\Lambda)\otimes W$
and in order to discuss characteristic identities it is important to emphasize that acting on $W$
the entries of $A$ are defined by
$$
A(e_\alpha\otimes w) = e_\beta\otimes A_{\beta\alpha}w \mbox{ (sum on $\beta$)}
$$
since then we may define matrix powers in the usual way:
$$
A^m(e_\alpha\otimes w) = e_\beta\otimes\left(A^m\right)_{\beta\alpha}w
$$
where 
$$
\left( A^m \right)_{\beta\alpha} = \sum_\gamma A_{\beta\gamma}\left( A^{m-1} \right)_{\gamma\alpha}
= \sum_\gamma \left(A^{m-1}\right)_{\beta\gamma}A_{\gamma\alpha}.
$$
\underline{{\em Remark}:} On the other hand if we have a representation of $A$ as
$$
A = \sum_{\alpha,\beta}e_{\alpha\beta}\otimes \hat{A}_{\alpha\beta},
$$ 
then
\begin{align*}
A(e_\gamma\otimes w) &= \sum_{\alpha,\beta}e_{\alpha\beta}e_\gamma\otimes
\hat{A}_{\alpha\beta}w(-1)^{([\alpha]+[\beta])[\gamma]}\\
&= \sum_\alpha e_\alpha\otimes \hat{A}_{\alpha\gamma}w(-1)^{[\gamma]([\alpha]+[\gamma])}\\
\Rightarrow \ \ A_{\alpha\gamma} &= (-1)^{[\gamma]([\alpha]+[\gamma])}
\hat{A}_{\alpha\gamma}.
\end{align*}
Thus if $A_{\alpha\beta}$ are the entries of $A$ then we have the representation
$$
A = e_{\alpha\beta}\otimes A_{\alpha\beta} (-1)^{[\beta]([\alpha]+[\beta])}.
$$

\begin{flushright}$\Box$\end{flushright}

The matrix $A$ acting on an irreducible module $V(\mu)$ determines an invariant on the
tensor product space:
\begin{align*}
A &= (q-q^{-1})^{-1} \pi_\Lambda\otimes \pi_\mu (I\otimes I - R^TR)\\
&= (q-q^{-1})^{-1} \pi_\Lambda\otimes \pi_\mu (I\otimes I -(v\otimes v)\Delta(v^{-1})).
\end{align*}
If $\lambda_1,\ldots,\lambda_k$ are the {\em distinct} weights in $V(\Lambda)$ then $A$ satisfies
the polynomial identity
\begin{equation}
\prod_{i=1}^k(A-a_i)=0,
\label{p23star}
\end{equation}
where
$$
a_i = \frac{1-q^{2\alpha_i}}{q-q^{-1}},
$$
with
\begin{align*}
\alpha_i &= \frac12\left[
(\mu+\lambda_i,\mu+\lambda_i+2\rho)-(\mu,\mu+2\rho)-(\Lambda,\Lambda+2\rho) \right]\\
&= \frac12\left[ (\lambda_i,\lambda_i+2(\mu+\rho)) - (\Lambda,\Lambda+2\rho) \right],
\end{align*}
which are the classical characteristic roots \cite{Gould1987}. This is the quantum analogue of the
characteristic identities considered by Green, Gould and others
\cite{Green1971,BraGre1971,Gould1978}. We now consider the
case $L= gl(m|n)$ which will be our focus for the remainder of the paper.


\subsection{$L$-operator for $U_q[gl(m|n)]$}\label{Section9} 

Let $\pi$ be the (undeformed) vector representation defined by
\begin{align*}
\pi(e_i) &= \pi(E_{i\ i+1}) = e_{i\ i+1},\\
\pi(f_i) &= \pi(E_{i+1\ i}) = e_{i+1\ i},\\
\pi(E_{ii}) &= e_{ii},
\end{align*}
with $e_{ij}$ and elementary matrix. Then we have the following $L$-operator (c.f. 
\cite{GI2019}):
$$
(\pi\otimes\id)R = q^{h_i\otimes h^i}\left\{ I\otimes I + (q-q^{-1})\sum_{i<j}e_{ji}\otimes
\hat{E}_{ij} \right\}
$$
where $\hat{E}_{ij}$ is given recursively by
$$
\hat{E}_{ij} = (-1)^{[k]}\left( \hat{E}_{ik}\hat{E}_{kj} -
q^{-(\varepsilon_k,\varepsilon_k)}\hat{E}_{kj}\hat{E}_{ik} \right),\ \ i<k<j
$$
with 
$$
\hat{E}_{i\ i+1} = (-1)^{[i]}q^{\frac12(\varepsilon_{i+1},\varepsilon_{i+1})} E_{i\ i+1}q^{h_i/2}.
$$
Above we have 
\begin{align*}
h_i = h_{\alpha_i} &= (-1)^{[i]}(E_{ii} - (-1)^{[\alpha_i]}E_{i+1\ i+1})\\
&= (-1)^{[i]}E_{ii} - (-1)^{[i+1]}E_{i+1\ i+1}\\
&=(i)E_{ii} - (i+1)E_{i+1\ i+1}
\end{align*}
where
\begin{equation}
(i) = (-1)^{[i]} = \left\{ \begin{array}{rl} 1 &, \mbox{ $i$ even}\\
-1 &, \mbox{ $i$ odd} \end{array}
\right.
=(\varepsilon_i,\varepsilon_i),
\label{p24star}
\end{equation}
so
\begin{align*}
h_i &= E_{ii} - E_{i+1\ i+1},\ \ 1\leq i< m,\\
h_m &= E_{mm} + E_{m+1\ m+1} = E_{mm} + E_{\mu=1\ \mu=1},\ \ (\alpha_m \mbox{ odd simple root}),\\
h_\mu &= E_{\mu+1\ \mu+1} - E_{\mu\mu},\ \ 1\leq \mu<n.
\end{align*}
We find it convenient to write 
$$
R = \sum_{i\leq j} e_{ji}\otimes \tilde{E}_{ij},
$$
where
$$
\tilde{E}_{ij} = 
\left\{ 
\begin{array}{rl}  
(q-q^{-1})E_{ij}', & i<j\\
q^{(j)E_{jj}}, & i=j
\end{array} 
\right.
$$
with 
$$
E_{ij}' = q^{(j)E_{jj}}\hat{E}_{ij},\ \ i<j.
$$
Then the recursion relation for the $\hat{E}_{ij}$ implies the following recursion
relation
$$
E_{ij}' = (-1)^{[k]}(E_{ik}'E_{kj}' - E_{kj}'E_{ik}')q^{-(k)E_{kk}}
$$
with
\begin{align*}
E_{i\ i+1}' &= q^{(i+1)E_{i+1\ i+1}} \hat{E}_{i\ i+1}\\
&= (-1)^{[i]}q^{\frac12((i)E_{ii} + (i+1)E_{i+1\ i+1} - (i))}E_{i\ i+1}.
\end{align*}
Finally, this suggests we set 
$$
E_{ij}' = (-1)^{[i]}q^{\frac12((i)E_{ii} + (j)E_{jj} - (i))}E_{ij}
$$
so that
$$
E_{i\ i+1}' = (-1)^{[i]}q^{\frac12((i)E_{ii} + (i+1)E_{i+1\ i+1} - (i))}E_{i\ i+1}
$$
as above. Then the recursion relation for the $E_{ij}'$ gives the following recursion relation for the
$E_{ij}$:
\begin{align}
E_{ij} &= E_{ik}E_{kj} - q^{-(k)}E_{kj}E_{ik} \label{p25star}\\
&= E_{ik}E_{kj} - q^{-(\varepsilon_k,\varepsilon_k)}E_{kj}E_{ik},\ \ i<k<j.\nonumber
\end{align} 
This gives the following ansatz for the $L$-operator:
\begin{equation}
(\pi\otimes\id)R = \sum_{i\leq j} e_{ji}\otimes \tilde{E}_{ij}
\label{p25starstar}
\end{equation}
with
$$
\tilde{E}_{ij} = \left\{ 
\begin{array}{rl}
(q-q^{-1})(-1)^{[i]}q^{\frac12((i)E_{ii} + (j)E_{jj} - (i))}E_{ij} &,\ i<j,\\
q^{(i)E_{ii}} &,\ i=j,
\end{array}
\right.
$$
with $E_{ij}$ given by the recursion relation (\ref{p25star}) and with
$(i)=(\varepsilon_i,\varepsilon_i)$ as in (\ref{p24star}). We refer to Appendix \ref{AppendixB} for a
discussion.


\subsection{Opposite $L$-operator}\label{Section10} 

We have
\begin{align*}
(\pi\otimes \id)R^T &= (\pi\otimes\id)R^\dagger\\
&= \sum_{a\leq b} e_{ba}^\dagger \otimes \tilde{E}_{ab}^\dagger (-1)^{[a]+[b]}\\
&= \sum_{a\leq b}e_{ab}\otimes \tilde{E}_{ba}
\end{align*}
where
$$
\tilde{E}_{ba} = \tilde{E}_{ab}^\dagger (-1)^{[a]+[b]},\ \ a\leq b
$$
with $\dagger$ the conjugation operation defined by 
$$
e_i^\dagger = f_i,\ \ h_i^\dagger = h_i.
$$
This gives immediately, for $i<j$, 
$$
\tilde{E}_{ji} = (q-q^{-1})(-1)^{[j]}q^{\frac12((i)E_{ii}+(j)E_{jj} - (j))}E_{ij}^\dagger
$$
where we now define $E_{ji} = E_{ij}^\dagger$, given recursively as follows:
\begin{align*}
E_{ji} &= E_{ij}^\dagger\\
&= (E_{ik}E_{kj} - q^{-(\varepsilon_k,\varepsilon_k)}E_{kj}E_{ik})^\dagger\\
&= E_{jk}E_{ki} - q^{-(\varepsilon_k,\varepsilon_k)}E_{ki}E_{jk},\ \ i<k<j.
\end{align*}
Thus we arrive at
\begin{equation}
(\pi\otimes\id)R^T = \sum_{i\leq j}e_{ij}\otimes \tilde{E}_{ji}
\label{p26star}
\end{equation}
with
\begin{equation}
\tilde{E}_{ji}
= \left\{  \begin{array}{rl}
(q-q^{-1})(-1)^{[j]}q^{\frac12((i)E_{ii}+(j)E_{jj}-(j))}E_{ji} &  , i<j,\\
q^{(i)E_{ii}} & , i=j,
\end{array}
\right.
\label{p27star}
\end{equation}
with $E_{ji}$ given recursively by
$$
E_{ji} = E_{jk}E_{ki} - q^{-(\varepsilon_k,\varepsilon_k)}E_{ki}E_{jk},\ \ i\lessgtr k\lessgtr j
$$
with
$$
E_{i\ i+1} = e_i,\ \ E_{i+1\ i}= f_i.
$$


\subsection{$A$-matrix}\label{Section11} 

The entries of the corresponding $A$-matrix,
$$
\hat{A} = (q-q^{-1})^{-1}(I\otimes I - (\pi\otimes \id)R^TR)
$$
follow from
\begin{align*}
(\pi\otimes \id)R^TR &= (\pi\otimes\id)R^T\cdot(\pi\otimes \id)R \\
&=\sum_{i\leq \ell}e_{i\ell}\otimes \tilde{E}_{\ell i}\cdot\sum_{j\leq k}e_{kj}\otimes
\tilde{E}_{jk}\\
&= \sum_{i\leq \ell}\sum_{j\leq k}(-1)^{([i]+[\ell])([j]+[k])}\delta_{k\ell}e_{ij}\otimes 
\tilde{E}_{\ell i}\tilde{E}_{jk}\\
&= \sum_{k\geq i\vee j} e_{ij}\otimes  \tilde{E}_{k i}\tilde{E}_{jk}(-1)^{([i]+[k])([j]+[k])},
\end{align*}
where $i\vee j = \mbox{max}(i,j)$. The entries $\hat{A}_{ij}$ of $\hat{A}$ are defined by
\begin{align*}
\hat{A} &= \sum_{i,j} e_{ij}\otimes \hat{A}_{ij} (-1)^{[j]([i]+[j])}\\
\Rightarrow \ \ \hat{A}_{ij} &= (q-q^{-1})^{-1}(-1)^{[j]([i]+[j])}
\left\{\delta_{ij} - \sum_{k\geq i\vee j} \tilde{E}_{ki}\tilde{E}_{jk} (-1)^{([i]+[k])([j]+[k])}
\right\}.
\end{align*}


\subsection{Transformation properties}\label{Section12} 

From the discussion above, the operators $\hat{A}_{ij}$ transform as a tensor operator of rank
$\pi^*\otimes \pi$, i.e. an adjoint tensor operator:
$$
a\hat{A}_{ij} = (-1)^{[a_{(2)}][i] + [a_{(3)}]([i]+[j])} \pi^*(a_{(1)})_{ki}\pi(a_{(2)})_{\ell
j}\hat{A}_{k\ell}a_{(3)}.
$$
If we set $i=j=m+n$ into the above, we obtain for $a\in U_q[gl(m|n-1)]$,
\begin{align*}
a\hat{A}_{m+n\ m+n} &= (-1)^{[a_{(2)}]} \pi^*(a_{(1)})_{k\ m+n}\pi(a_{(2)})_{\ell \ m+n}
\hat{A}_{k\ell}a_{(3)}\\
&= \hat{A}_{m+n\ m+n}a
\end{align*}
so $\hat{A}_{m+n\ m+n}$ is an invariant of $U_q[gl(m|n-1)]$.

Now set $j=m+n$ and consider the transformation properties of the last column, $\hat{A}_{i\ m+n}$,
$1\leq i\leq m+n$ under $U_q[gl(m|n-1)]$: for $a\in U_q[gl(m|n-1)]$ we have
\begin{align*}
a\hat{A}_{i\ m+n} &= (-1)^{[a_{(2)}][i]+[a_{(3)}]([i]+1)}\pi^*(a_{(1)})_{ki}
\pi(a_{(2)})_{\ell\ m+n} \hat{A}_{k\ell}a_{(3)}
&= (-1)^{[a_{(2)}]([i]+1)} \pi^*(a_{(1)})_{ki}\hat{A}_{k\ m+n}a_{(2)}.
\end{align*}
Therefore 
$$
\phi_i\equiv (-1)^{[i]} \hat{A}_{i\ m+n}
$$
transforms as an {\em odd dual vector} operator:
$$
a \phi_i = (1)^{[a]+[a_{(2)}][i]} \pi^*(a_{(1)})_{ki}\phi_k a_{(2)}.
$$
Similarly for the last row $\hat{A}_{m+n\ i}$ we have
\begin{align*}
a\hat{A}_{m+n\ i} &= (-1)^{[a_{(2)}]+[a_{(3)}](1+[i])} \pi^*(a_{(1)})_{k\ m+n}
\pi(a_{(2)})_{\ell i}\hat{A}_{k\ell} a_{(3)}\\
&= (-1)^{[a_{(1)}]+[a_{(2)}](1+[i])}  \pi(a_{(1)})_{\ell i}\hat{A}_{m+n\ \ell} a_{(2)}\\
&= (-1)^{[a]+[a_{(2)}][i]} \pi(a_{(1)})_{\ell i} \hat{A}_{m+n\ \ell} a_{(2)}
\end{align*}
which is the transformation law for an {\em odd vector} operator.


\subsection{Polynomial identities}\label{Section13} 

In this case the identities (\ref{p23star}) reduce to 
$$
\prod_{r=1}^{m+n}(\hat{A} - \overline{a}_r)=0
$$
where 
$$
\overline{a}_i = \frac{1-q^{2\overline{\alpha}_i}}{q-q^{-1}}
$$
with
$$
\overline{\alpha}_i = \frac12[(\varepsilon_i,\varepsilon_i + 2(\Lambda+\rho)) -
(\varepsilon_1,\varepsilon_1+2\rho)].
$$

Using the form of $\rho$ from (\ref{p18star}), we obtain the classical adjoint roots
\begin{align}
\overline{\alpha}_i &= \Lambda_i + 1-i,\ \ 1\leq i\leq m, \label{p30stara}\\
\overline{\alpha}_\mu &= \mu-m-1 - \Lambda_\mu,\ \ 1\leq\mu\leq n. \label{p30starb}
\end{align} 
\underline{{\em Note}:} It is worth noting that if $\overline{\rho}$ is the graded half
sum of the positive roots of $gl(m|n-1)$ then
\begin{align*}
\rho-\overline{\rho} &= \frac12\sum_{\mu=1}^n\delta_\mu - \frac12\sum_{i=1}^m
\varepsilon_i + \frac12(m-n)\delta_n\\
&= \frac12(m-n)\delta_n - \frac12\sum_a(-1)^{[a]} \varepsilon_a.
\end{align*}
Therefore for $i,j<n$ we have
$$
(\rho-\overline{\rho},\varepsilon_i-\varepsilon_j) =
\frac12(-1)^{[j]}(\varepsilon_j,\varepsilon_j) -
\frac12(-1)^{[i]}(\varepsilon_i,\varepsilon_i) = 0.
$$

Following the classical procedure our aim is to apply the above identities to the
construction of projection operators and the evaluation of invariants. However, as for
$U_q[gl(n)]$ the first $m+n-1$ rows and columns of the matrix $\hat{A}$ do not reduce to
the corresponding matrix of $U_q[gl(m|n-1)]$. This leads us to consider the $R$-matrix
$\tilde{R} = (R^T)^{-1}$.


\subsection{$L$-operators for $\tilde{R}$}\label{Section14} 

Using
\begin{align*}
\tilde{R} &= (R^T)^{-1} = ((S\otimes\id)R)^T = (\id\otimes S)R^T,\\
\tilde{R}^T &= R^{-1} = (\id\otimes S^{-1})R,
\end{align*}
we obtain the following $L$-operators:
\begin{align*}
(\pi\otimes\id)\tilde{R} &= (\pi\otimes\id)(\id\otimes S)R^T\\
& \stackrel{(\ref{p26star})}{=} \sum_{i\leq j} e_{ij}\otimes S(\tilde{E}_{ji}),\\
(\pi\otimes\id)\tilde{R^T} &= (\pi\otimes\id)(\id\otimes S^{-1})R\\
& \stackrel{(\ref{p25starstar})}{=} \sum_{i\leq j} e_{ji}\otimes S^{-1}(\tilde{E}_{ij}),
\end{align*}
with $\tilde{E}_{ij}$ as in (\ref{p27star}).


\subsection{Matrix $\tilde{A}$}\label{Section15} 

The corresponding characteristic matrix is
$$
\tilde{A} = (q-q^{-1})^{-1}\left( I - (\pi\otimes\id)\tilde{R}^T\tilde{R} \right)
$$
where 
$$
(\pi\otimes\id)\tilde{R}^T\tilde{R} = \sum_{k\leq i\wedge
j}(-1)^{([k]+[i])([k]+[j])} e_{ij}\otimes
S^{-1}(\tilde{E}_{ki})S(\tilde{E}_{jk}),
$$
with $i\wedge j=\mbox{min}(i,j)$.
The entries of this matrix defined by
$$
\tilde{A} = \sum_{i,j} e_{ij}\otimes \tilde{A}_{ij}(-1)^{[j]([i]+[j])}
$$
are given by 
$$
\tilde{A}_{ij} = (q-q^{-1})^{-1} (-1)^{[j]([i]+[j])}\left\{ \delta_{ij} - \sum_{k\leq
i\wedge j} (-1)^{([k]+[i])([k]+[j])} S^{-1} (\tilde{E}_{ki}) S(\tilde{E}_{jk}) \right\}.
$$
The entries of this matrix transform in the same way as the matrix $\hat{A}$ but now has
the desirable feature that the block containing the first $m+n-1$ rows and columns reduce
to the characteristic matrix of $U_q[gl(m|n-1)]$, as in the $U_q[gl(n)]$ and classical
case. Hence we work with the matrix $\tilde{A}$ rather than $\hat{A}$.

Thus for the last column of $\tilde{A}$, the operators
$$
\phi_i \equiv (-1)^{[i]}\tilde{A}_{i\ m+n},\ \ i<m+n
$$
transform as an odd dual vector operator, while the entries of the last row,
$$
\psi_i \equiv \tilde{A}_{m+n\ i},
$$
transform as an odd vector operator, and the $(m+n,m+n)$ entry $\tilde{A}_{m+n\ m+n}$ is a
$U_q[gl(m|n-1)]$ scalar.

The characteristic identity satisfied by $\tilde{A}$ is
$$
\prod_{r=1}^{m+n} (\tilde{A} - \tilde{a}_r) = 0,
$$
where (c.f. \cite{GI2019})
$$
\tilde{a}_r = \frac{1-q^{-2\bar{\alpha}_r}}{q-q^{-1}},
$$
with $\bar{\alpha}_r$ as in equations (\ref{p30stara}) and (\ref{p30starb}). Here we make
the simplifying assumption that these roots are distinct on $V(\Lambda)$, and the set of
$\Lambda\in D_+$ for which this occurs is Zariski dense in $D_+$. With this assumption, if
$V$ is the vector module, we have the irreducible $U_q[gl(m|n)]$-module direct sum
$$
V\otimes V(\Lambda) = \bigoplus_{r=1}^{m+n} V(\Lambda+\varepsilon_r).
$$
The projections
$$
\tilde{P}_r = \prod_{i\neq r}^{m+n} \left( \frac{\tilde{A} - \tilde{a}_i}{\tilde{a}_r -
\tilde{a}_i} \right)
$$
project onto the submodule $V(\Lambda+\varepsilon_r)$.


\subsection{Dual vector $L$-operator}\label{Section16} 

We have
\begin{align*}
	(\pi^*\otimes\id)R &= (\pi^T\otimes\id)(S\otimes\id)R =
	(\pi^T\otimes\id)(\id\otimes S^{-1})R\\
	(\pi^*\otimes\id)R^T &= (\pi^T\otimes\id)(S\otimes\id)R^T =
	(\pi^T\otimes\id)(\id\otimes S^{-1})R^T,
\end{align*}
where, as before, $T$ is super-transpose, so
$$
e_{ij}^T = (-1)^{[i]([i]+[j])}e_{ji}.
$$
We thus obtain
\begin{align}
(\pi^*\otimes\id)R &\stackrel{(\ref{p25starstar})}{=} \sum_{i\leq j}(-1)^{[j]([i]+[j])}e_{ij}\otimes
S^{-1}(\tilde{E}_{ij}), \label{p34stara}\\
(\pi^*\otimes\id)R^T &\stackrel{(\ref{p26star})}{=} \sum_{i\leq j}(-1)^{[i]([i]+[j])}e_{ji}\otimes
S^{-1}(\tilde{E}_{ji}). \label{p34starb}
\end{align}
Thus we arrive at the characteristic matrix 
\begin{align*}
A &= (q-q^{-1})^{-1}[I - (\pi^*\otimes \id)R^TR] \\
&=(q-q^{-1})^{-1}\left\{ I - \sum_{k\leq i\wedge j}(-1)^{[j]([i]+[j])}
e_{ij} \otimes S^{-1}(\tilde{E}_{ik}) \tilde{S}^{-1}(\tilde{E}_{kj})\right\} 
\end{align*}
with entries 
\begin{align}
A_{ij} = (q-q^{-1})^{-1}\left\{ \delta_{ij} - \sum_{k\leq i\wedge j}
S^{-1}(\tilde{E}_{ik})S^{-1}(\tilde{E}_{kj}) \right\}
\label{p34starstar}
\end{align}

\underline{{\em Remark}:} From the point of view of reduced Wigner
coefficients, to be investigated in a later section, it is important to note
that although the vector representation $\pi$ is (type I) unitary, the dual
$\pi^*$ is not (type II) unitary - moreover $\pi^*$ is {\em not} undeformed.
The corresponding undeformed unitary irreducible representation is defined by 
$$
\overline{\pi}(E_{ij}) = -E_{ji},\ \ |i-j|\leq 1
$$
and is given by 
$$
\overline{\pi}(a) = \pi^T(\gamma(a)),\ \ a\in U_q[gl(m|n)],
$$
where
$$
\gamma(a) = q^{h_\rho}S(a)q^{-h_\rho}
$$
is the {\em principal anti-automorphism}:
\begin{align*}
\gamma(x) &= -x,\ \ x=e_i, f_i, h_i \\
	\Rightarrow \ \ \overline{\pi}(a) &= \pi^T(q^{h_\rho}S(a)q^{-h_\rho})\\
	&= \pi^T(S(q^{h_\rho}aq^{-h_\rho})) = \pi^*(q^{h_\rho}aq^{-h_\rho}).
\end{align*}

\begin{flushright}$\Box$\end{flushright}

The matrix $A$ satisfies the following polynomial identity on $V(\Lambda)$:
$$
\prod_{r=1}^{m+n}(A-a_r)=0,\ \ a_r = \frac{1-q^{-2\alpha_r}}{q-q^{-1}},
$$
where
\begin{align*}
	\alpha_r &= -\frac12\left[ (2(\Lambda+\rho-\varepsilon_r,-\varepsilon)
	- (2\rho-\delta_n,-\delta_n)\right] \\
	\Rightarrow \ \ \alpha_i &= \Lambda_i+m-n-i,\ \ 1\leq i\leq m \mbox{
		(even roots),} \\
		\alpha_\mu &= \mu-n-\Lambda_\mu,\ \ 1\leq \mu\leq n  \mbox{
			(odd roots)}.
\end{align*}
In this case the entries of the matrix $A$ transform as a tensor operator of
rank $\pi^{**}\otimes\pi^*.$

\underline{{\em Note}:} Using the following property of super transposition 
$$
(A^T)^T = (-1)^{[A]}A
$$
it follows that
\begin{align*}
	\pi^{**}(a) &= (-1)^{[a]}\pi(S^2(a))\\
	&= (-1)^{[a]} \pi(q^{-2h_\rho}aq^{2h_\rho}).
\end{align*}
Thus we have the transformation law
\begin{align}
	aA_{ij} &=
	(-1)^{[a_{(2)}][i]+[a_{(3)}]([i]+[j])}\pi^{**}(a_{(1)})_{ki}\pi^*(a_{(2)})_{\ell
	j} A_{k\ell}a_{(3)} \nonumber \\
	&= (-1)^{[a_{(1)}]+[a_{(2)}][i] + [a_{(3)}]([i]+[j])}
	q^{2(\rho,\varepsilon_i-\varepsilon_k)}
	\pi(a_{(1)})_{ki}\pi^*(a_{(2)})_{\ell j} A_{k\ell}a_{(3)}.
	\label{p36star}  
\end{align}

\noindent
\underline{{\em Note}:} Thus $\hat{\hat{A}}_{ij} \equiv
q^{-2(\rho,\varepsilon_i)}A_{ij}$ transforms as
$$
a \hat{\hat{A}}_{ij} = (-1)^{[a_{(1)}]+[a_{(2)}][i]+[a_{(3)}]([i]+[j])}
\pi(a_{(1)})_{ki} \pi^*(a_{(2)})_{\ell j}  \hat{\hat{A}}_{k\ell} a_{(3)}.
$$


\subsection{Partitioning of $A$}\label{Section17} 

As for the matrix $\tilde{A}$, $A$ may be partitioned as
$$
A = \left( \begin{array}{ccc:c} 
&&& \\ 
\ & A_0 & & * \\ 
&&& \\ 
\hdashline
&&& \\
 & * & & A_{m+n\ m+n} \end{array} \right)
$$
with $A_0$ the matrix of $U_q[gl(m|n-1)]$. Setting $i=j=m+n$ into
(\ref{p36star}) gives, for $a\in U_q[gl(m|n-1)]$,
\begin{align*}
	aA_{m+n\ m+n} &= (-1)^{[a_{(1)}]+[a_{(2)}]}
	q^{2(\rho,\varepsilon_{m+n}-\varepsilon_k)} \pi(a_{(1)})_{k\
	m+n}\pi^*(a_{(2)})_{\ell\ m+n} A_{k\ell} a_{(3)}\\
	&= A_{m+n\ m+n}a,
\end{align*}
so that $A_{m+n\ m+n}$ is a $U_q[gl(m|n-1)]$ scalar.

\underline{{\em Remark}:} It is worth noting the following matrix elements for $a\in
U_q[gl(m|n-1)]$ which we will repeatedly use:
$$
\pi(a)_{k\ m+n} = \pi(a)_{m+n\ k} = \varepsilon(a) \delta_{m+n\ k}.
$$
Similarly
$$
\pi^*(a)_{k\ m+n} = \pi^*(a)_{m+n\ k} = \varepsilon(a) \delta_{m+n\ k},
$$
since, for example,
\begin{align*}
\pi^*(a)_{k\ m+n} &= \langle a e_{m+n}^*,e_k \rangle \\
&= \langle e_{m+n}^*, S(a) e_k\rangle (-1)^{[a]} \\
&= \delta_{k\ m+n}\varepsilon(a) (-1)^{[a]} = \delta_{k\ m+n} \varepsilon(a).
\end{align*}

\begin{flushright}$\Box$\end{flushright}

Next set $j=m+n$ to give, for $a\in U_q[gl(m|n-1)]$,
\begin{align}
	a A_{i\ m+n} &= (-1)^{[a_{(1)}]+[a_{(2)}][i]+[a_{(3)}]([i]+1)}
	q^{2(\rho,\varepsilon_i-\varepsilon_k)} \pi(a_{(1)})_{ki}
	\pi^*(a_{(2)})_{\ell \ m+n} A_{k\ell} a_{(3)}\nonumber \\
	&= (-1)^{[a]} q^{2(\rho,\varepsilon_i-\varepsilon_k)}
	(-1)^{[a_{(2)}][i]} \pi(a_{(1)})_{ki} A_{k\ m+n} a_{(3)}\nonumber\\
	\Rightarrow \ \ \psi^i &= q^{-2(\rho,\varepsilon_i)} A_{i\ m+n}
	\label{p37star}
\end{align}
transforms as an odd vector operator, i.e.
$$
a\psi^i = (-1)^{[a]} (-1)^{[a_{(2)}][i]} \pi(a_{(1)})_{ki} \psi^k a_{(3)}.
$$
Finally, set $i=m+n$ to give, for $a\in U_q[gl(m|n-1)]$
\begin{align*}
aA_{m+n\ j} &= (-1)^{[a_{(1)}] + [a_{(2)}] + [a_{(3)}](1+[j])} q^{2(\rho,\varepsilon_{m+n}
- \varepsilon_k)} \pi(a_{(1)})_{k\ m+n} \pi^*(a_{(2)})_{\ell j} A_{k\ell} a_{(3)}\\
&= (-1)^{[a]+[a_{(2)}][j]} \pi^*(a_{(1)})_{\ell j} A_{m+n\ j} a_{(2)} 
\end{align*}
so the $A_{m+n\ i}$ ($i<m+n$) transforms as an odd dual vector operator.

The matrix $A$ gives an operator on the tensor product space $V^*\otimes V(\Lambda)$.
Below we assume $\Lambda$ is such that the classical characteristic roots are all
distinct, so we have a direct sum of irreducible $U_q[gl(m|n)]$-modules
$$
V^*\otimes V(\Lambda) = \bigoplus_{r=1}^{m+n} V(\Lambda-\varepsilon_r).
$$
The corresponding projection operators
$$
P_r = \prod_{\ell\neq r}^{m+n} \left( \frac{A-a_\ell}{a_r-a_\ell} \right)
$$
project onto the module $V(\Lambda-\varepsilon_r)\subset V^*\otimes V(\Lambda).$

\noindent
\underline{{\em Remark}:} We may make this assumption clear by considering the polynomial
function on $H^*$ defined by
$$
p(\lambda) = \prod_{i<j}^{m+n} (\alpha_i-\alpha_j)(\bar{\alpha}_i-\bar{\alpha}_j).
$$ 
Then
$$
K = \ker p = \{ \lambda\in H^*|p(\lambda) = 0\}
$$
is closed in the Zariski topology. Therefore its complement $\tilde{K}$ is open, and thus
dense in the Zariski topology. We may assume $\Lambda\in D_+\cap\tilde{K}$ is Zariski
dense in $D_+$. Any polynomial function derived on $\tilde{K}$ can be extended by
continuity to all $H^*$.


\subsection{Summary of vector operators}\label{Section18} 

Below we shall utilise the partitioning of the matrices $A$ and $\tilde{A}$,
together with the characteristic identities, to evaluate certain invariants.
In the table that follows, we summarise the transformation properties of the $A_{ij}$,
$\tilde{A}_{ij}$, where $i$ or $j$ is $m+n$, under $U_q[gl(m|n-1)]$:

\begin{center}
\begin{tabular}{|c|c|}
  \hline
	& \\
	$A_{m+n\ m+n}$, $\tilde{A}_{m+n\ m+n}$ & $U_q[gl(m|n-1)]$ scalars \\
	& \\
	\hline
	& \\
	$\tilde{\phi}_i=(-1)^{[i]}\tilde{A}_{i\ m+n}$ & odd dual vector
	operator \\ 
	& \\
	\hline
	& \\
	$\tilde{\psi}_i=\tilde{A}_{m+n\ i}$ & odd vector
	operator \\ 
	& \\
	\hline
	& \\
	$\psi_i=q^{-2(\rho,\varepsilon_i)}A_{i\ m+n}$ & odd vector
	operator \\ 
	& \\
	\hline
	& \\
	$\phi_i=A_{m+n\ i}$ & odd dual vector
	operator \\ 
	& \\
	\hline
\end{tabular}
\end{center}

\noindent
\underline{{\em Note}}: If $\overline{\rho}$ is the half-sum of positive roots
for $gl(m|n-1)$ we have seen that
$$
\rho = \overline{\rho}+\Delta,\ \ \Delta = -\frac12\sum_{a=1}^{m+n-1}
(-1)^{[a]}\varepsilon_a+\frac12(m-n)\delta_{n}. 
$$
Therefore for $i<m+n$ we have
\begin{align*}
(\Delta,\varepsilon_i) &= -\frac12(-1)^{[i]}(\varepsilon_i,\varepsilon_i) = -\frac12\\
\Rightarrow \ \ (2\rho,\varepsilon_i) &= (2\overline{\rho},\varepsilon_i)-1\\
\Rightarrow \ \ \psi_i &= q^{-2(\rho,\varepsilon_i)}A_{i\ m+n}\\ 
	&= q^{1-2(\overline{\rho},\varepsilon_i)}A_{i\ m+n}.
\end{align*}
Therefore if we replace $\rho$ above with $\overline{\rho}$ we still obtain an
odd vector operator under $U_q[gl(m,n-1)]$.


\subsection{Shift components}\label{Section19} 

For quantum supergroups we have the following result (c.f. Lemma \ref{lem14}).
\begin{lemma}
\label{lemma40}
Let $V$ be a finite dimensional irreducible $U_q(L)$-module with homogeneous basis $\{e_\alpha\}$.
Then (sum on $\alpha$)
\begin{align*}
v_0 &= e_\alpha\otimes e_\alpha^*,\\
\overline{v}_0 &= (-1)^{[\alpha]} e_\alpha^*\otimes q^{2h_\rho}e_\alpha
\end{align*}
span the identity submodules of $V\otimes V^*$, $V^*\otimes V$, respectively.
\end{lemma} 

Now suppose $\{\phi_i\}$ is a homogeneous dual vector operator of $U_q[gl(m|n)]$ on
$V(\Lambda)$ that determines an intertwining operator on $V^*\otimes V(\Lambda)$, i.e.
$$
\phi(e_i^*\otimes v) = \phi_iv,\ \ v\in V(\Lambda).
$$ 
Using the identity resolution
$$
I = \sum_{r=1}^{m+n}P_r
$$
implies
$$
\phi = \sum_{r=1}^{m+n}\phi[r],
$$
where the shift components $\phi[r]$ are given by
\begin{align*}
\phi[r]_iv &= \phi[r](e_i^*\otimes v)\\
&= \phi\left( P_r(e_i^*\otimes v) \right) \\
&= \phi(e_j^*\otimes (P_r)_{ji}v) \\
&= \phi_j (P_r)_{ji}v\\
\Rightarrow \ \ \phi[r]_i &= \phi_j(P_r)_{ji} \mbox{ (sum on $j$).}
\end{align*}

Clearly
$$
\phi[r]V^*\otimes V(\Lambda) \subseteq V(\Lambda-\varepsilon_r)
$$
so for $v\in V(\Lambda)$,
$$
\phi[r]_i v\in V(\Lambda-\varepsilon_r)
$$
which implies $\phi[r]$ is a pure shift dual vector operator.

On the other hand, we may consider the projection $\tilde{P}_r$ of $V\otimes V(\Lambda)$
onto $V(\Lambda +\varepsilon_r)$. With $v_0$ as in the lemma above we note that
$$
(\id \otimes \phi)(\mathbb{C}v_0\otimes V(\Lambda)) \cong V(\Lambda)
$$
and for $v\in V(\Lambda)$
\begin{align*}
(\id \otimes \phi)(v_0\otimes v) &= (\id\otimes \phi)(e_i\otimes e_i^*\otimes v)\mbox{
(sum on $i$)}\\
&= (-1)^{[\phi][i]} e_i\otimes \phi_iv.
\end{align*}
Then we note that
$$
\tilde{P}_r(\id\otimes \phi[\ell])(v_0\otimes v) \subseteq
V(\Lambda-\varepsilon_\ell+\varepsilon_r)\cap V(\Lambda)
$$
$$
\Rightarrow \ \ \tilde{P}_r(\id\otimes\phi[\ell])(v_0\otimes v) = (0),\ \ r\neq \ell.
$$
Hence we have
\begin{align*}
(-1)^{[\phi][i]}e_i\otimes \phi[r]_iv 
&= e_i\otimes \phi[r](e_i^*\otimes v)(-1)^{[\phi][i]}\\
&= \tilde{P}_r(\id\otimes \phi)(v_0\otimes v)\\
&= \tilde{P}_r(e_i\otimes \phi(e_i^*\otimes v))(-1)^{[\phi][i]} \\
&= e_j\otimes (\tilde{P}_r)_{ji}\phi_iv (-1)^{[\phi][i]} \\
\Rightarrow \ \ (-1)^{[j][\phi]}\phi[r]_j &= (\tilde{P}_r)_{ji} \phi_i (-1)^{[\phi][i]} \mbox{ (sum
on $i$)}
\end{align*}

In this sense the projection $\tilde{P}_r$ projects out shift components from the left.
Thus, for even dual vector operators, we have
$$
\phi[r]_i = (\tilde{P}_r)_{ij}\phi_j
$$
as in the usual case, and for the odd case we have
\begin{equation}
(-1)^{[i]} \phi[r]_i = (\tilde{P}_r)_{ij} (-1)^{[j]}\phi_j. 
\label{p42star}
\end{equation}

On the other hand, given a homogeneous vector operator $\psi^i$, $\tilde{P}_r$ projects
out the shift components of $\psi^i$ from the {\em right}:
$$
\psi[r]_i = \psi_j (\tilde{P}_r)_{ji} \mbox{ (sum on $j$)}
$$
and for $v\in V(\Lambda)$
$$
\psi[r]_iv = \psi \tilde{P}_r(e_i\otimes v)\in V(\Lambda+\varepsilon_r),
$$
i.e. $\psi[r]$ increases the highest weight by $\varepsilon_r$.

As above, we consider $\bar{v}_0$ as in the lemma and note the isomorphism
$$
(\id\otimes\psi)(\mathbb{C}\bar{v}_0\otimes V(\Lambda)) \cong V(\Lambda).
$$
As above, we now note that
$$
P_r(\id\otimes \psi[\ell])(\bar{v}_0\otimes v)\in V(\Lambda+\varepsilon_\ell -
\varepsilon_r)\cap V(\Lambda)
$$
$$
\Rightarrow \ \ P_r(\id\otimes \psi[\ell])(\bar{v}_0\otimes V(\Lambda)) = (0),\ \ \ell\neq
r.
$$
Hence we now obtain for $v\in V(\Lambda)$
\begin{align*}
(-1)^{[i]([i]+[\psi])} q^{(2\rho, \varepsilon_i)}e_i^*\otimes \psi[r]_iv
&= (-1)^{[i]([i]+[\psi])} e_i^*\otimes \psi[r] (q^{2h_\rho}e_i\otimes v)\\
&= (\id\otimes \psi[r])(\bar{v}_0\otimes v)\\
&= P_r(\id\otimes \psi)(\bar{v}_0\otimes v)\\
&= (-1)^{[i]([i]+[\psi])} P_r(e_i^*\otimes \psi(q^{2h_\rho}e_i\otimes v) )\\
&= (-1)^{[i]([i]+[\psi])}q^{2(\rho, \varepsilon_i)} P_r(e_i^*\otimes \psi_iv)\\
&= (-1)^{[i]([i]+[\psi])}q^{2(\rho, \varepsilon_i)} e_j^*\otimes (P_r)_{ji}\psi_iv
\end{align*}
\begin{equation}
\Rightarrow \ \ (-1)^{[i]([i]+[\psi])}q^{2(\rho, \varepsilon_i)} \psi[r]_i
= (-1)^{[j]([j]+[\psi])} (P_r)_{ij} q^{2(\rho, \varepsilon_j)} \psi_j.
\label{p43star}
\end{equation}
Therefore in this sense, $P_r$ projects out shift components from the left.

\noindent
\underline{{\em Remark}:} The above suggests we introduce the matrix
$$
\overline{A}_{ij} = (-1)^{[i]([i]+[\psi])}q^{-2(\rho, \varepsilon_i)} A_{ij}
(-1)^{[j]([j]+[\psi])} q^{2(\rho, \varepsilon_j)}
$$
which satisfies the same polynomial identity as $A$ and whose projections $\overline{P}_r$
project out shift components of a vector operator $(\psi)$ from the left.

The above matrix $\overline{A}$ seems to be related to the characteristic matrix
associated with the undeformed dual vector representation $\overline{\pi}$. We consider
this in more detail later.


\subsection{Invariants}\label{Section20} 

Let us go back to the general case and consider an irreducible $U_q(L)$-module $V(\Lambda)$ with
homogeneous basis $\{e_\alpha\}$ and $\pi_\Lambda$ the representation afforded by $V(\Lambda)$.
Recall that the entries of the corresponding characteristic matrix 
$$
A = (q-q^{-1})^{-1}\left\{ I - (\pi_\Lambda\otimes \id)R^TR \right\}
$$
are defined via
$$
A(e_\alpha\otimes w) = e_\beta\otimes A_{\beta\alpha}w,\ \ w\in W
$$
acting on a given $U_q(L)$-module $W$. Moreover from Proposition \ref{prop16} we have an
intertwining operator $\varphi$ on $V(\Lambda)^*\otimes V(\Lambda)\otimes W$ defined by
$$
\varphi(e_\alpha^*\otimes e_\beta\otimes w) = A_{\alpha\beta}w.
$$
Then we have seen that
$$
a\varphi(e_\alpha^*\otimes e_\beta\otimes w) = \varphi\left( (\id\otimes\Delta)\Delta(a)
(e_\alpha^*\otimes e_\beta\otimes w) \right)
$$
or
$$
aA_{\alpha\beta} w = \varphi\left( (\id\otimes\Delta)\Delta(a)(e_\alpha^*\otimes
e_\beta\otimes w) \right).
$$
Now consider (c.f. Lemma \ref{lemma40})
\begin{align*}
\bar{v}_0 &= (-1)^{[\alpha]} e_\alpha^* \otimes q^{2h_\rho} e_\alpha\\
&=(-1)^{[\alpha]}q^{2(\rho,\lambda_\alpha)}  e_\alpha^* \otimes e_\alpha,
\end{align*}
where it is understood that we sum on $\alpha$ and $\lambda_\alpha$ is the weight of basis
vector $e_\alpha$. Then $\bar{v}_0$ spans the identity submodule of $V(\Lambda)^*\otimes
V(\Lambda)$ and by direct calculation
\begin{align*}
\varphi(\bar{v}_0\otimes w) &= (-1)^{[\alpha]} q^{2(\rho,\lambda_\alpha)}
\varphi(e_\alpha^*\otimes e_\alpha\otimes w)\\
&= C w, 
\end{align*}
where 
$$
C = (-1)^{[\alpha]} q^{2(\rho,\lambda_\alpha)} A_{\alpha\alpha} \mbox{ (sum on $\alpha$).}
$$
Further we observe that for $a\in U_q(L)$
\begin{align*}
aCw &= a\varphi(\bar{v}_0\otimes w) \\
&= \varphi( (\id\otimes\Delta)\Delta(a)(\bar{v}_0\otimes w) ) \\
&= \varphi( \Delta(a_{(1)})\bar{v}_0\otimes a_{(2)} w) \\
&= \varphi(\bar{v}_0\otimes aw)\\
&= Ca w
\end{align*}
or 
$$
aC = Ca,\ \ \forall a\in U_q(L).
$$
Thus
$$
C = \mbox{str}_\Lambda\left( q^{2h_\rho}A \right)
$$
is an invariant. This fact applies also to any power of the matrix $A$, so we arrive at
the family of invariants
$$
C_m^\Lambda = \mbox{str}_\Lambda\left( q^{2h_\rho}A^m \right).
$$
In particular we obtain the following invariants arising from the $U_q[gl(m|n)]$ matrices
$A$ and $\tilde{A}$:
\begin{align*}
C_m &= \mbox{str}\left( \pi^*(q^{2h_\rho})A^m \right)\\
&= (-1)^{[i]} q^{-2(\rho,\varepsilon_i)} A_{ii}^m \mbox{ (sum on $i$),} \\
\tilde{C}_m &= \mbox{str}\left( \pi(q^{2h_\rho}) \tilde{A}^m \right) \\
&= (-1)^{[i]} q^{2(\rho,\varepsilon_i)} \tilde{A}_{ii} \mbox{ (sum on $i$).}
\end{align*}
In particular we have the first order invariants
\begin{align*}
C_1 &= (-1)^{[i]} q^{-2(\rho,\varepsilon_i)} A_{ii},\\
\tilde{C}_1 &= (-1)^{[i]} q^{2(\rho,\varepsilon_i)} \tilde{A}_{ii} \mbox{ (sum on $i$),}
\end{align*}
where
\begin{align*}
A_{ii} &= (q-q^{-1})^{-1} \left\{ I - \sum_{k\leq i} S^{-1}(\tilde{E}_{ik})
S^{-1}(\tilde{E}_{ki}) \right\},\\
\tilde{A}_{ii} &= (q-q^{-1})^{-1} \left\{ I - \sum_{k\leq i} (-1)^{[i]+[k]}
S^{-1}(\tilde{E}_{ki}) S(\tilde{E}_{ik}) \right\}.
\end{align*} 

Note that if $v_+^\Lambda$ is a maximal weight vector of weight $\Lambda$, then
$$
A_{ii} v_+^\Lambda = (q-q^{-1})^{-1} \left\{ I - S^{-1}(\tilde{E}_{ii})
S^{-1}(\tilde{E}_{ii}) \right\} v_+^\Lambda.
$$
Using 
$$
\tilde{E}_{ii} = q^{(i) E_{ii} } \ \ \Rightarrow \ \ S(\tilde{E}_{ii}) = q^{-(i)E_{ii}}
$$
we have
\begin{align*}
A_{ii} v_+^\Lambda &= (q-q^{-1})^{-1} \left\{ I - q^{-2(i)E_{ii}} \right\} v_+^\Lambda\\
&= (q-q^{-1})^{-1} \left\{ 1 - q^{-2(i)\Lambda_i} \right\} v_+^\Lambda\\
&= (q-q^{-1})^{-1} \left\{ 1 - q^{-2(\Lambda,\varepsilon_i)} \right\} v_+^\Lambda.
\end{align*}
Therefore the eigenvalue of $C_1$ on $V(\Lambda)$ is given by
\begin{align*}
\chi_\Lambda(C_1) &= \sum_{i=1}^{m+n} (-1)^{[i]} q^{-2(\rho,\varepsilon_i)}
\frac{1-q^{-2(\Lambda,\varepsilon_i)}}{q-q^{-1}} \\
&= \sum_{i=1}^{m+n} (-1)^{[i]} q^{-(\Lambda + 2\rho,\varepsilon_i)} 
[(\Lambda,\varepsilon_i) ]_q 
\end{align*}
where we have used the $q$-number notation
$$
[x]_q = \frac{q^x-q^{-x}}{q-q^{-1}}.
$$
Similarly let $v_-^\Lambda$ be the minimal weight vector of the Kac module $K(\Lambda)$,
given by the weight vector of weight
$$
\Lambda_- = \Lambda_-^0 - 2\rho_1
$$
with $\Lambda_-^0 = \tau(\Lambda)$, the minimal weight of the $L_0 =
gl(m)\oplus gl(n)$ irreducible module $V_0(\Lambda)$. Then we have
\begin{align*}
\tilde{A}_{ii} v_-^\Lambda
&= (q-q^{-1})^{-1} \left\{ I - S^{-1}(\tilde{E}_{ii}) S(\tilde{E}_{ii}) \right\}
v_-^\Lambda\\
&= (q-q^{-1})^{-1} \left\{ I - q^{-2(i)E_{ii}} \right\} v_-^\Lambda \\
&= (q-q^{-1})^{-1} \left\{ 1 - q^{-2(\Lambda_-^0 - 2\rho_1,\varepsilon_i)}  \right\}
v_-^\Lambda. 
\end{align*}
Therefore on $V(\Lambda)$ the invariant $\tilde{C}_1$ takes the eigenvalue
\begin{align*}
\chi_\Lambda (\tilde{C}_1) &= \sum_i (-1)^{[i]} q^{2(\rho,\varepsilon_i)} 
\frac{1-q^{-2(\Lambda_-^0-2\rho_1,\varepsilon_i)}}{q-q^{-1}} \\
&= \sum_i (-1)^{[i]} q^{2(\rho,\tau(\varepsilon_i))}
\frac{1-q^{-2(\Lambda_-^0-2\rho_1,\tau(\varepsilon_i))}}{q-q^{-1}}. 
\end{align*}
Now
\begin{align*}
(\rho,\tau(\varepsilon_i)) &= (\rho_0-\rho_1,\tau(\varepsilon_i)) =
-(\rho_0+\rho_1,\varepsilon_i) \\
(\Lambda_-^0-2\rho_1,\tau(\varepsilon_i)) &= (\Lambda - 2\rho_1,\varepsilon_i),
\end{align*}
which leads to
\begin{align*}
\chi_\Lambda (\tilde{C}_1)
&= \sum_i (-1)^{[i]} q^{-2(\rho_0+\rho_1,\varepsilon_i)}
\frac{1-q^{-2(\Lambda-2\rho_1,\varepsilon_i)}}{q-q^{-1}}\\
&= \sum_i (-1)^{[i]} q^{-2(\rho_0+\rho_1,\varepsilon_i) - (\Lambda -
2\rho_1,\varepsilon_i)} [(\Lambda-2\rho_1,\varepsilon_i)]_q \\
&= \sum_i (-1)^{[i]} q^{-(\Lambda-2\rho_0,\varepsilon_i)}
[(\Lambda-2\rho_1,\varepsilon_i)]_q.
\end{align*}


\subsection{The matrix $\overline{A}$}\label{Section21} 

Here we determine the $L$-operator and associated characteristic matrix arising from the undeformed
dual vector representation $\overline{\pi}$. We recall that
$$
\overline{\pi}(a) = \pi^*(q^{h_\rho}aq^{-h_\rho}). 
$$
We have from equations (\ref{p34stara}) and (\ref{p34starb}) 
\begin{align*}
(\overline{\pi}\otimes\id)R   &= \sum_{i\leq j}
(-1)^{[j]([i]+[j])}q^{(\rho,\varepsilon_j-\varepsilon_i)} e_{ij}\otimes
S^{-1}(\tilde{E}_{ij}) \\
(\overline{\pi}\otimes\id)R^T &=\sum_{i\leq j}
(-1)^{[i]([i]+[j])}q^{(\rho,\varepsilon_i-\varepsilon_j)} e_{ji}\otimes
S^{-1}(E_{ji}) 
\end{align*}
which gives the characteristic matrix
\begin{align*}
\overline{A} &= (q-q^{-1})^{-1}\left\{ I - (\overline{\pi}\otimes\id)RR^T \right\} \\
&= (q-q^{-1})^{-1}\left\{ I - \sum_{k\leq i\wedge j}
q^{(\rho,\varepsilon_j-\varepsilon_i)} (-1)^{[j]([i]+[j])} e_{ij}\otimes
S^{-1}(\tilde{E}_{ik})S^{-1}(\tilde{E}_{kj}) \right\}
\end{align*}
with entries 
$$
\overline{A}_{ij} =  (q-q^{-1})^{-1}q^{(\rho,\varepsilon_j-\varepsilon_i)}
\left\{ I - \sum_{k\leq i\wedge j}
(-1)^{[j]([i]+[j])} e_{ij}\otimes
S^{-1}(\tilde{E}_{ik})S^{-1}(\tilde{E}_{kj}) \right\}.
$$
This is related to the entries (\ref{p34starstar}) of the matrix $A$ by
\begin{equation}
\overline{A}_{ij} = q^{(\rho,\varepsilon_j-\varepsilon_i)} A_{ij} =
q^{-(\rho,\varepsilon_i)} A_{ij} q^{(\rho,\varepsilon_j)}
\label{p49star}
\end{equation}
so satisfies the same polynomial identity. These operators now transform as a tensor
operator of rank $\overline{\pi}^*\otimes \overline{\pi},$ i.e.
$$
a\overline{A}_{ij} = (-1)^{[a_{(2)}][i]+[a_{(3)}]([i]+[j])}
\overline{\pi}^*(a_{(1)})_{ki}\overline{\pi}(a_{(2)})_{\ell j} \overline{A}_{k\ell}
a_{(3)},
$$
where we note that
\begin{align*}
\overline{\pi}(a) &= \pi^*(q^{h_\rho} a q^{-h_\rho}) = (-1)^{[a]} \pi( q^{-h_\rho} a
q^{h_\rho} ).
\end{align*}
We thus arrive at the transformation law 
$$
a\overline{A}_{ij} = (-1)^{[a_{(1)}] + [a_{(2)}][i]+[a_{(3)}]([i]+[j]) } 
q^{(\rho, \varepsilon_i+\varepsilon_j - \varepsilon_k - \varepsilon_\ell )} 
\pi(a_{(1)})_{ki} \pi^*(a_{(2)})_{\ell j} \overline{A}_{k\ell} a_{(3)},
$$
which is consistent, via equation (\ref{p49star}), with the result of equation
(\ref{p36star}). 

\noindent
\underline{{\em Note}:} If $\overline{\rho}$ is the graded half sum of positive roots for
$gl(m|n-1)$, we have seen that
$$
(\rho-\overline{\rho},\varepsilon_i-\varepsilon_j) = 0 \mbox{ for } 1\leq i< j<n,
$$
so actually $\overline{A}$ forms blocks in the required way under $U_q[gl(m|n-1)].$
Setting $i=j=m+n$ into the above gives the following transformation law for $a\in
U_q[gl(m|n-1)]$:
$$
a\overline{A}_{m+n\ m+n} = \overline{A}_{m+n\ m+n}a.
$$
Similarly setting $i=m+n$ gives
\begin{align*}
a\overline{A}_{m+n\ j} &= (-1)^{[a_{(1)}]+[a_{(2)}]+[a_{(3)}](1+[j])}
q^{(\rho,\varepsilon_j-\varepsilon_\ell)} \pi(a_{(1)})_{m+n\ m+n} \pi^*(a_{(2)})_{\ell j}
\overline{A}_{m+n\ \ell} a_{(3)} \\
&= (-1)^{[a]+[a_{(2)}][j]} q^{(\rho,\varepsilon_j-\varepsilon_\ell)} \pi^*(a_{(1)})_{\ell j}
\overline{A}_{m+n\ \ell} a_{(2)}\\
&=  (-1)^{[a]+[a_{(2)}][j]} \overline{\pi}(a_{(1)})_{\ell j} \overline{A}_{m+n\ \ell}
a_{(2)},\ \ a\in U_q[gl(m|n-1)],
\end{align*}
so $\overline{A}_{m+n\ j}$ is a dual vector operator of type $\overline{\pi}$, i.e. a
pseudo vector operator.

Similarly setting $j=m+n$ gives, for $a\in U_q[gl(m|n-1)]$, 
\begin{align*}
a\overline{A}_{i\ m+n} &= (-1)^{[a_{(2)}][i]+[a_{(3)}]([i]+1)+[a_{(1)}]}
q^{(\rho,\varepsilon_i+\varepsilon_{m+n} - \varepsilon_k-\varepsilon_\ell)}
\pi(a_{(1)})_{ki}\pi^*(a_{(2)})_{\ell\ m+n} \overline{A}_{k\ell} a_{(3)}\\
&= (-1)^{[a_{(1)}]+[a_{(2)}]([i]+1)} q^{(\rho,\varepsilon_i-\varepsilon_k)}
\pi(a_{(1)})_{ki}\overline{A}_{k\ m+n} a_{(2)}\\
&= (-1)^{[a]+[a_{(2)}][i]}  q^{(\rho,\varepsilon_i-\varepsilon_k)}
\pi(a_{(1)})_{ki}\overline{A}_{k\ m+n} a_{(2)}.
\end{align*}
This implies that
$$
\overline{\psi}_i = q^{-(\rho,\varepsilon_i)} \overline{A}_{i\ m+n}
$$
transforms as an {\em odd} vector operator.

\noindent
\underline{{\em Note}:} $\overline{\phi}=q^{-(\rho,\varepsilon_j)} \overline{A}_{m+n\ j}$
transforms as a dual vector operator.


\section{Reduced Wigner coefficients}\label{Section4} 

\subsection{Index sets}

We first recall the $gl(m|n)\downarrow gl(m|n-1)$ branching condition. In a
$U_q[gl(m|n)]$-module $V(\Lambda)$ the allowed $gl(m|n-1)$ irreducible
submodules have highest weights $\Lambda_0$ such that
$$
\Lambda_\mu\geq {\Lambda_0}_\mu\geq \Lambda_{\mu+1} \mbox{ (Gelfand-Tsetlin condition)}
$$
and 
$$
{\Lambda_0}_i = \Lambda_i \mbox{ or } \Lambda_i-1,\ \ \mbox{ i.e. }
\Lambda_i\geq {\Lambda_0}_i\geq \Lambda_i-1.
$$
Thus it is suggestive that we introduce the index sets
\begin{align*}
	I_0(\Lambda,\Lambda_0) &= \{ i=1,\ldots,m|{\Lambda_0}_i = \Lambda_i-1
	\},\\
	\overline{I}_0(\Lambda,\Lambda_0) &= \{
		i=1,\ldots,m|{\Lambda_0}_i=\Lambda_i\}
\end{align*}
so $I_0\cup \overline{I}_0 = \{i=1,\ldots,m\}$.

These index sets are important since, following the classical case, the operators
$$
\tilde{\omega}_r = \left( \tilde{P}_r \right)_{m+n\ m+n},\ \ 
\omega_r = \left( P_r \right)_{m+n\ m+n}
$$
we have seen, are $U_q[gl(m|n-1)]$ invariants, whose eigenvalues determine squares of
reduced Wigner coefficients:
$$
\tilde{\omega}_k = 
\left| 
\left\langle \begin{array}{c} \Lambda+\varepsilon_k\\ \Lambda_0 \end{array} \right| 
\left| \begin{array}{c} \Lambda\\ \Lambda_0 \end{array} \right\rangle
\right|^2, \ \ 
\omega_k = 
\left| 
\left\langle \begin{array}{c} \Lambda-\varepsilon_k\\ \Lambda_0 \end{array} \right| 
\left| \begin{array}{c} \Lambda\\ \Lambda_0 \end{array} \right\rangle
\right|^2
$$
so we must have
$$
\tilde{\omega}_k = 0 \mbox{ if } k\in I_0,\ \ 
\omega_k = 0  \mbox{ if } k\in \overline{I}_0.
$$
Moreover, if $\psi$ is a vector operator, $\phi$ a dual vector operator of
$U_q[gl(m|n-1)]$ acting on $V(\Lambda)$, with shift components $\psi[r]$, $\phi[r]$
respectively, we must have
$$
\psi[r]=0 \mbox{ if } r\in \overline{I}_0,\ \ 
\phi[r]=0 \mbox{ if } r\in I_0.
$$

\subsection{Evaluation of invariants}\label{Section23} 

Here we use the notation employed in \cite{GI2019}. Using the $U_q[gl(m|n)]$
identity we have
\begin{align}
a_k (P_k)_{i\ m+n} 
&= (AP_k)_{i\ m+n}\nonumber\\
&= (A_0)_{ij}(P_k)_{j\ m+n} + A^i_{\ m+n}\omega_k\nonumber\\
\Rightarrow \ \ A^i_{\ m+n}\omega_k &= (a_k-A_0)_{ij} (P_k)_{j\ m+n}. \label{p52star} 
\end{align}
Now recall from (\ref{p37star}) that 
$$
\psi^i = q^{-2(\rho,\varepsilon_i)}A_{i\ m+n}
$$
is an odd vector operator and hence, from the $[\psi]=1$ case of equation (\ref{p43star}),
its shift components are given by  
\begin{align*}
q^{2(\rho,\varepsilon_i)} \psi[r]_i
&= \left( {P_0}_r\right)_{ij} q^{2(\rho,\varepsilon_j)} \psi_j\\
&= \left( {P_0}_r\right)_{ij} A_{j\ m+n}.
\end{align*}
Hence multiplying equation (\ref{p52star}) on the left by $({P_0}_r)_{ij}$ gives
\begin{align*}
q^{2(\rho,\varepsilon_i)} \psi[r]_i \omega_k 
&= (a_k-{a_0}_r)\left( {P_0}_r\right)_{ij} \left( P_k\right)_{j\ m+n} \\
\Rightarrow \ \  q^{2(\rho,\varepsilon_i)} (a_k-{a_0}_r)^{-1} \psi[r]_i \omega_k 
&= \left( {P_0}_r\right)_{ij} \left( P_k\right)_{j\ m+n}.
\end{align*}
Therefore summing on $r$ gives
$$
\sum_r (a_k-{a_0}_r)^{-1} \psi[r]_i \omega_k = q^{-2(\rho,\varepsilon_i)}\left(
P_k\right)_{i\ m+n}.
$$
Now we note that
\begin{align*}
{\alpha_0}_r \psi[r]_i
&= \psi[r]_i ({\alpha_0}_r + (-1)^{[r]})\\
&= \psi[r]_i ({\alpha_0}_r + (r) ),\ \ \left((r)=(-1)^{[r]}\right)
\end{align*}
where ${\alpha_0}_r$, ${P_0}_r$, $A_0$ are the classical characteristic roots, projections
and $A$-matrix for the subalgebra $U_q[gl(m|n-1)]$.

Using
$$
{a_0}_r = \frac{1-q^{-2{\alpha_0}_r}}{q-q^{-1}}
$$
we obtain
$$ 
{a_0}_r\psi[r]_i = \psi[r]_i \left( q^{-2(r)} {a_0}_r + (r)q^{-(r)} \right).
$$
Therefore we have 
\begin{equation}
\sum_r \psi[r]_i \left( a_k - q^{-2(r)} {a_0}_r - (r)q^{-(r)} \right)^{-1} \omega_k
= q^{-2(\rho,\varepsilon_i)} \left( P_k \right)_{i\ m+n}.
\label{p53star}
\end{equation}
Therefore summing over $k$, we obtain
$$
\sum_{r,k} \psi[r]_i\left( a_k - q^{-2(r)} {a_0}_r - (r)q^{-(r)} \right)^{-1} \omega_k =0.
$$

Now observe that $\psi[r]$, $\omega_k$ vanish when $r,k\in \overline{I}_0$. Thus if we let
$I_1$ denote the odd index set $I_1=\{ \mu=1,\ldots,n-1 \}$ and set $\tilde{I}_1 = I_1\cup
\{m+n\}$, then the linear independence of the $\psi[r]_i$ gives the set of equations
$$
\sum_{k\in I_0\cup\tilde{I}_1} \left( a_k - q^{-2(r)} {a_0}_r - (r)q^{-(r)} \right)^{-1}
\omega_k =0,
r\in I_0\cup I_1.
$$
This gives $|I_0|+n-1$ equations in $|I_0|+n$ unknowns $\omega_k$. These equations,
together with 
$$
\sum_{k\in I_0\cup \tilde{I}_1} \omega_k = \sum_{k\in I_0\cup \tilde{I}_1} \left(
P_k \right)_{m+n\ m+n} = 1
$$
uniquely determine the $\omega_k$.
Following our previous work we readily obtain
$$
\omega_k = \prod_{r\in I_0\cup I_1}\left( a_k - q^{-2(r)} {a_0}_r - (r)q^{-(r)} \right) 
\prod_{\ell\neq k}^{I_0\cup\tilde{I}_1} (a_k-a_\ell)^{-1},\ \ k\in I_0\cup\tilde{I}_1,
$$
where $\omega_k=0$ for $k\in \overline{I}_0.$ Thus the explicit form of the eigenvalues
depends on the index set $I_0$, which in turn depends on the $U_q[gl(m|n)]$ highest
weight $\Lambda$ and $U_q[gl(m|n-1)]$ highest weight $\Lambda_0$.

Similarly applying the identity for the matrix $\tilde{A}$, we have, in obvious subalgebra
notation as above,
$$
\tilde{a}_k\left( \tilde{P}_k \right)_{i\ m+n} = \left( \tilde{A} \tilde{P}_k \right)_{i\ m+n}
= \left( \tilde{A}_0 \right)_{ij}\left( \tilde{P}_k \right)_{j\ m+n}+\tilde{A}_{i\
m+n}\tilde{\omega}_k, 
$$
where we recall that 
$$
\phi_i = (-1)^{[i]}\tilde{A}_{i\ m+n}
$$
is an odd vector operator. From equation (\ref{p42star}), the shift components are given
by
$$
(-1)^{[i]} \phi[r]_i = \left( \tilde{P}_{0r} \right)_{ij} (-1)^{[j]} \phi_j = \left(
\tilde{P}_{0r} \right)_{ij} \tilde{A}_{j\ m+n}.
$$ 
Therefore multiplying the above equation on the left by $\tilde{P}_{0r}$ gives
$$
\left( \tilde{a}_k - \tilde{a}_{0r} \right) \left( \tilde{P}_{0r} \right)_{ij} 
\left( \tilde{P}_k \right)_{j\ m+n} = (-1)^{[i]} \phi[r]_i \tilde{\omega}_k
$$
$$
\Rightarrow \ \ \left( \tilde{P}_{0r} \right)_{ij}  \left( \tilde{P}_k \right)_{j\ m+n}
= (-1)^{[i]} \left( \tilde{a}_k - \tilde{a}_{0r} \right)^{-1} \phi[r]_i \tilde{\omega}_k.
$$
Therefore summing on $r$ gives
$$
\sum_{r=1}^{m+n} \left( \tilde{a}_k - \tilde{a}_{0r} \right)^{-1} \phi[r]_i \tilde{\omega}_k
= \left( \tilde{P}_k \right)_{i\ m+n}(-1)^{[i]}.
$$
Now using 
$$
\tilde{a}_{0r} = \frac{1-q^{-2\bar{\alpha}_{0r}}}{q-q^{-1}}
$$
together with the shift
$$
\bar{\alpha}_{0r} \phi[r]_i = \phi[r]_i \left( \bar{\alpha}_{0r} - (r) \right)
$$
we obtain
$$
\tilde{a}_{0r} \phi[r]_i = \phi[r]_i \left( q^{2(r)} \tilde{a}_{0r} - (r) q^{(r)} \right)
$$
\begin{equation}
\Rightarrow \ \ \sum_{r=1}^{m+n} \phi[r]_i\left( \tilde{a}_k -q^{2(r)} \tilde{a}_{0r} +
(r) q^{(r)}\right)^{-1} \tilde{\omega}_k = (-1)^{[i]} \left( \tilde{P}_k \right)_{i\ m+n}.
\label{p55star}
\end{equation}
Therefore summing over $k$ we arrive at the equations
$$
\sum_{r,k} \phi[r]_i\left( \tilde{a}_k -q^{2(r)} \tilde{a}_{0r} +
(r) q^{(r)}\right)^{-1} \tilde{\omega}_k = 0.
$$

Now we observe for $r,k\in I_0$ that $\phi[r]_i$, $\tilde{\omega}_k$ vanish. For
$r\in\overline{I}_0$ the $\phi[r]_i$ are linearly independent so we obtain the equations
$$
\sum_{k\in\overline{I}_0\cup\tilde{I}_1} \left( \tilde{a}_k -q^{2(r)} \tilde{a}_{0r} +
(r) q^{(r)}\right)^{-1} \tilde{\omega}_k = 0,\ \ r\in \overline{I}_0\cup I_1,
$$
which, together with 
$$
\sum_{k\in\overline{I}_0\cup\tilde{I}_1} \tilde{\omega}_k = \sum_{k=1}^{m+n}
\tilde{\omega}_k = 1,
$$
uniquely determines the $\tilde{\omega}_k$. We thus obtain
$$
\tilde{\omega}_k = \prod_{r\in\overline{I}_0\cup I_1} 
\left( \tilde{a}_k -q^{2(r)} \tilde{a}_{0r} + (r) q^{(r)}\right)
\prod_{\ell\neq k}^{\overline{I}_0\cup\tilde{I}_1} \left( \tilde{a}_k - \tilde{a}_\ell
\right)^{-1}.
$$


\subsection{Reduced matrix elements: Summary and notation}\label{Section24} 

We have the odd dual vector operators
\begin{align*}
\tilde{\phi}_i &= (-1)^{[i]} \tilde{A}_{i\ m+n},\\
\phi_i &= A_{m+n\ i},
\end{align*}
where $P_r$ projects out shift components from the right and $\tilde{P}_r(-1)^{[i]}$ from the left.
We also have the odd vector operators
\begin{align*}
\tilde{\psi}_i &= \tilde{A}_{m+n\ i},\\
\psi_i &= q^{-2(\rho,\varepsilon_i)} A_{m+n\ i},
\end{align*}
where $\tilde{P}_r$ projects out shift components from the right and
$q^{-2(\rho,\varepsilon_j)}\left( P_r \right)_{ji}q^{2(\rho,\varepsilon_i)}$ from the left.

Now from equation (\ref{p53star}) we have 
\begin{equation}
\left( P_k \right)_{i\ m+n} = q^{2(\rho,\varepsilon_i)} \sum_{r\in I_0\cup I_1} \psi[r]_i
\left( a_k-q^{-2(r)} {a_0}_r - (r)q^{-(r)} \right)^{-1} \omega_k, \ \ k\in I_0\cup
\tilde{I}_1.
\label{p56star}
\end{equation}
Now we attempt to invert this equation by looking for the unique solution $\gamma_{rk}$,
$r\in I_0\cup I_1$, $k\in I_0\cup \tilde{I}_1$ to the equations
\begin{align}
\sum_{k\in I_0\cup \tilde{I}_1} \gamma_{rk} \left( a_k - q^{-2(\ell)}{a_0}_\ell -
(\ell)q^{-(\ell)} \right)^{-1}\omega_k &= \delta_{r\ell}
\label{p56starstara}
\\ 
\sum_{k\in I_0\cup \tilde{I}_1} \gamma_{rk} \omega_k &= 0,\ \ r,\ell\in I_0\cup I_1. 
\label{p56starstarb}
\end{align}

For each $r\in I_0\cup I_1,$ this yields $|I_0|+n$ equations in $|I_0|+n$ unknowns
$\gamma_{rk}$, $k\in I_0\cup \tilde{I}_1$ with unique solution
$$
\gamma_{rk} = \gamma_r\left( a_k - q^{-2(r)} {a_0}_r - (r) q^{-(r)} \right)^{-1}
$$
where
$$
\gamma_r = (-1)^{|I_0|+n-1} \prod_{k\in I_0\cup \tilde{I}_1} \left( a_k - q^{-2(r)} {a_0}_r - (r) q^{-(r)} \right)
\prod_{\ell\neq r}^{I_0\cup I_1} \left( q^{-2(r)} {a_0}_r - q^{-2(\ell)}{a_0}_\ell +
(r)q^{-(r)} - (\ell) q^{-(\ell)} \right)^{-1}.
$$
These invariants essentially determine the {\em $q$-length} of the vector operator
$\psi_i$.
Indeed, multiplying equation (\ref{p56star}) by $\gamma_{rk}$ and summing on $k$ gives
\begin{align*}
\sum_{k}^{I_0\cup \tilde{I}_1} \left(P_k\right)_{i\ m+n} \gamma_{rk}
&= q^{2(\rho,\varepsilon_i)} \sum_\ell^{I_0\cup I_1}\sum_k^{I_0\cup \tilde{I}_1}
\psi[\ell]_i \left( a_k - q^{-2(\ell)} {a_0}_\ell - (\ell) q^{-(\ell)} \right)^{-1}
\omega_k\gamma_{rk} \\
&\stackrel{(\ref{p56starstara})}{=} q^{2(\rho,\varepsilon_i)} \psi[r]_i.
\end{align*}
Therefore multiplication on the left by $\phi_i = A_{m+n\ i}$ gives
\begin{align*}
\phi[r]_i q^{2(\rho,\varepsilon_i)} 
&= \phi_i q^{2(\rho,\varepsilon_i)} \psi[r]_i \mbox{ (sum on $i$)}
\\
&= \sum_k^{I_0\cup\tilde{I}_1} A_{m+n\ i}\left( P_k \right)_{i\ m+n} \gamma_{rk}
\\
&= \sum_k^{I_0\cup\tilde{I}_1} \left( a_k - A_{m+n\ m+n} \right) \omega_k \gamma_{rk}
\\
&\stackrel{(\ref{p56starstarb})}{=} \sum_k a_k\omega_k\gamma_{rk}
\\
&= \gamma_r\sum_k a_k\omega_k\left( a_k - q^{-2(r)} {a_0}_r - (r) q^{-(r)} \right)^{-1}
\\
&= \gamma_r\sum_k \left\{ \omega_k + \left( q^{-2(r)} {a_0}_r + (r)q^{-(r)} \right)
\omega_k \left( a_k - q^{-2(r)} {a_0}_r - (r) q^{-(r)} \right)^{-1} \right\}
\\
&\stackrel{(\ref{p53star})}{=} \gamma_r.
\end{align*}
It follows that $q^{2(\rho,\varepsilon_i)}\psi[r]_i(\gamma_r)^{-1}\phi[r]_j$ is a zero
shift tensor transforming as $\pi^{**}\otimes \pi^*$ so we must have
$$
q^{2(\rho,\varepsilon_i)} \psi[r]_i (\gamma_r)^{-1} \phi[r]_j = \left( P_r\right)_{ij}.
$$
Now we note the root shift
$$
{a_0}_r \phi[r]_i = \phi[r]_i \left( q^{2(r)}{a_0}_r - (r) q^{(r)} \right),
$$
so that we may write
$$
\left( P_r\right)_{ij} = q^{2(\rho,\varepsilon_i)} \psi[r]_i \phi[r]_j \mu_r^{-1},
$$
where 
\begin{equation}
\mu_r = (-1)^{|i_0|+n-1} \prod_{k\in I_0\cup \tilde{I}_1} (a_k - {a_0}_r)
\prod_{\ell\neq r}^{I_0\cup I_1} \left( {a_0}_r - q^{-2(\ell)} {a_0}_\ell - (\ell)
q^{-(\ell)} \right)^{-1}.
\label{p58star}
\end{equation}
Thus we may write
\begin{equation}
q^{2(\rho,\varepsilon_i)} \psi[r]_i\phi[r]_j = \mu_r \left( P_r\right)_{ij}.
\label{p58starstar}
\end{equation}
The invariants $\mu_r$ determine the squared reduced matrix elements of the dual vector
operator $\phi[r]_i$.

Similarly, we have from equation (\ref{p55star})
$$
\sum_{r\in \overline{I}_0\cup I_1} \phi[r]_i \left( \tilde{a}_k - q^{2(r)}\tilde{a}_{0r} +
(r)q^{(r)} \right)^{-1} \tilde{\omega}_k = (-1)^{[i]} \left( \tilde{P}_k \right)_{i\ m+n}
$$
for $k\in \overline{I}_0\cup\tilde{I}_1.$ Now we consider the solution
$\tilde{\gamma}_{rk}$ to
\begin{align}
\sum_{k\in \overline{I}_0\cup\tilde{I}_1} \tilde{\gamma}_{rk} \tilde{\omega}_k \left(
\tilde{a}_k - q^{2(\ell)} \tilde{a}_{0\ell} + (\ell) q^{(\ell)} \right) &= \delta_{r\ell}, 
\label{p59stara} \\
\sum_{k\in \overline{I}_0\cup\tilde{I}_1} \tilde{\gamma}_{rk} \tilde{\omega}_k &= 0.
\label{p59starb}
\end{align} 
For each $r\in \overline{I}_0\cup I_1$, this gives $|\overline{I}_0|+n$ equations in
$|\overline{I}_0|+n$ unknowns $\tilde{\gamma}_{rk}$, $k\in\overline{I}_0\cup\tilde{I}_1$.
We have the unique solution
$$
\tilde{\gamma}_{rk} = \left( \tilde{a}_k - q^{-2(r)} \tilde{a}_{0r} + (r)q^{(r)} \right)^{-1} \tilde{\gamma}_r
$$
where
$$
\tilde{\gamma}_r = (-1)^{|\overline{I}_0|+n-1}
\frac{\prod_k^{\overline{I}_0\cup\tilde{I}_1}\left( \tilde{a}_k - q^{2(r)} \tilde{a}_{0r} +
(r) q^{(r)} \right)}{ \prod_{\ell\neq r}^{\overline{I}_0\cup I_1} \left(
q^{2(r)}\tilde{a}_{0r} - q^{2(\ell)} \tilde{a}_{0\ell} + (\ell) q^{(\ell)} - (r)q^{(r)}
\right) }.
$$
Now we multiply equation (\ref{p55star}) on the right by $\tilde{\gamma}_{rk}$ to give
\begin{align*}
(-1)^{[i]}\sum_{k\in \overline{I}_0\cup \tilde{I}_1} \left( \tilde{P}_k \right)_{i\ m+n}
\tilde{\gamma}_{rk}
&= \sum_{\ell\in\overline{I}_0\cup I_1} \phi[\ell]_i\left( \tilde{a}_k - q^{2(\ell)}
\tilde{a}_{0\ell} + (\ell) q^{(\ell)} \right)^{-1} \tilde{\omega}_k \tilde{\gamma}_{rk}\\
&= \phi[r]_i.
\end{align*}

At this point we change notation slightly and instead write the dual vector operator
$\phi_i$ as $\tilde{\phi}_i$ to be consistent with the expressions given at the start of
Section \ref{Section24}, i.e.
$$
\tilde{\phi}_i = (-1)^{[i]} \tilde{A}_{i\ m+n},\ \ i<m+n.
$$ 
Then the equation above becomes
$$
\tilde{\phi}[r]_i = (-1)^{[i]} \sum_{k\in \overline{I}_0 \cup \tilde{I}_1} \left(
\tilde{P}_k \right)_{i\ m+n} \tilde{\gamma}_{rk}.
$$
Therefore multiplying on the left by $\tilde{\psi}_i= \tilde{A}_{m+n\ i}$ and summing on $i$
gives
\begin{align*}
\tilde{\psi}[r]_i\tilde{\phi}[r]_i(-1)^{[i]} 
&= \tilde{\psi}_i (-1)^{[i]} \tilde{\phi}[r]_i \\
&= \sum_{k\in \overline{I}_0\cup \tilde{I}_1} \tilde{A}_{m+n\ i} \left( \tilde{P}_k
\right)_{i\ m+n} \tilde{\gamma}_{rk} \\
&= \sum_{k\in \overline{I}_0\cup \tilde{I}_1} \left( -\tilde{A}_{m+n\ m+n} + \tilde{a}_k \right)
\tilde{\omega}_k \tilde{\gamma}_{rk} \\
&\stackrel{(\ref{p59stara}),(\ref{p59starb})}{=} \tilde{\gamma}.
\end{align*} 
Thus in this case we may write
$$
\tilde{\phi}[r]_i \left( \tilde{\gamma}_r \right)^{-1} \tilde{\psi}[r]_j (-1)^{[i]} =
\left( \tilde{P}_r \right)_{ij}
$$
which may be rearranged to give
\begin{equation}
(-1)^{[i]} \tilde{\phi}[r]_i \tilde{\psi}[r]_j = \tilde{\mu}_r\left( \tilde{P}_r
\right)_{ij},
\label{p60star}
\end{equation}
where
$$
\tilde{\mu}_r(\Lambda,\Lambda_0) = \tilde{\gamma}_r(\Lambda,\Lambda_0+\varepsilon_r)
$$
$$
\Rightarrow \ \ \tilde{\mu}_r = (-1)^{|\overline{I}_0|+n-1} \frac{ \prod_k^{\overline{I}_0\cup\tilde{I}_1}
\left( \tilde{a}_k - \tilde{a}_{0r} \right) }{ \prod_{\ell\neq r}^{\overline{I}_0\cup I_1}
\left( \tilde{a}_{0r} - q^{2(\ell)} \tilde{a}_{0\ell} + (\ell) q^{(\ell)} \right) }.
$$
This formula determines the squared reduced matrix elements of the vector operator
$\tilde{\psi}_i = \tilde{A}_{m+n\ i}.$


\subsection{Squared reduced Wigner coefficients and other invariants}\label{Section25} 

Now from the characteristic identity we have
$$
\left( P_k \right)_{i\ m+n}A_{m+n\ j} = \left( P_k \right)_{i\ell}\left(a_k - A_0\right)_{\ell j}
$$
or
$$
\left( P_k \right)_{i\ m+n}\phi_j = \left( P_k \right)_{i\ell}\left(a_k - A_0\right)_{\ell j}.
$$
Therefore we can multiply on the right by $P_{0r}$ to give
$$
\left( P_k \right)_{i\ m+n} \phi[r]_j = \left( P_k \right)_{i\ell} \left( {P_0}_r
\right)_{\ell j} \left(a_k - {a_0}_r \right) 
$$
or
$$
\left( P_k {P_0}_r \right)_{ij} = \left( P_k \right)_{i\ m+n} \phi[r]_j \left(
a_k - a_{0r} \right)^{-1}.
$$
For a non-zero contribution on the left and right hand side we require $k\in
I_0\cup\tilde{I}_1$ and $r\in I_0\cup I_1.$

Now multiply the above on the left by ${P_0}_r$ to give
\begin{align*}
\left( {P_0}_r P_k {P_0}_r \right)_{ij} 
&=  \left( {P_0}_r\right)_{i\ell} \left(P_k \right)_{\ell\ m+n} \phi[r]_j
(a_k - a_{0r})^{-1} \\
&\stackrel{(\ref{p56star})}{=} 
q^{2(\rho,\varepsilon_i)} \psi[r]_i ( a_k - q^{-2(r)}a_{0r} - (r) q^{-(r)} )^{-1} 
\omega_k\phi[r]_j (a_k-a_{0r})^{-1}\\
&= \left( P_r\right)_{ij} \psi[r]_i q^{2(\rho,\varepsilon_i)} \phi[r]_j 
( a_k - q^{-2(r)}a_{0r} - (r) q^{-(r)} )^{-1} (a_k-a_{0r})^{-1} \omega_k
\end{align*}
since $( a_k - q^{-2(r)}a_{0r} - (r) q^{-(r)} )^{-1} \omega_k$ is independent of $a_{0r}$
and so commutes with $\phi[r]$. Finally, using equation (\ref{p58starstar}) we arrive at
$$
{P_0}_r P_k {P_0}_r = \omega_{kr} {P_0}_r,\ \ k\in I_0\cup\tilde{I}_1,\ r\in I_0\cup I_1,
$$ 
where
$$
\omega_{kr} = \omega_k\mu_r( a_k - q^{-2(r)}a_{0r} - (r) q^{-(r)} )^{-1}
(a_k-a_{0r})^{-1}.
$$
\underline{{\em Note}}: For $r\notin I_0\cup I_1$, we obtain zero for both sides of the
equation. If $r\in I_0\cup I_1$, but $k\notin I_0\cup \tilde{I}_1$ we also obtain zero
for both sides, but in that case we must have $\omega_{kr}=0$. Therefore we have a
non-zero contribution to the left and right hand sides only when $k\in
I_0\cup\tilde{I}_1$ and $r\in I_0\cup I_1$, as noted above.
\begin{flushright}$\Box$\end{flushright}

The $\omega_{kr}$ are of interest since they determine squared reduced Wigner coefficients:
$$
\chi_{(\Lambda,\Lambda_0)}(\omega_{kr}) = 
\left| 
\left\langle 
\left.
\begin{array}{c} \Lambda-\varepsilon_k\\ \Lambda_0-\varepsilon_{0r} \end{array}
\right|
\left|
\begin{array}{c} \overline{\varepsilon}_1\\ \overline{\varepsilon}_{1} \end{array}
;
\begin{array}{c} \Lambda\\ \Lambda_0 \end{array}
\right.
\right\rangle
\right|^2.
$$
Similarly for the matrix $\tilde{A}$ we have
$$
\left( \tilde{P}_k \right)_{i\ m+n} \tilde{A}_{m+n\ j} = \left( \tilde{P}_k
\right)_{i\ell}\left( \tilde{a}_k - \tilde{A}_0 \right)_{\ell j}.
$$
Therefore multiplying on the right by $\tilde{P}_{0r}$ gives
$$
\left( \tilde{P}_k \right)_{i\ m+n} \tilde{\psi}[r]_j = \left( \tilde{P}_k \right)_{i\ell}
\left( \tilde{P}_{0r} \right)_{\ell j} \left( \tilde{a}_k - \tilde{a}_{0r} \right),
$$
where we have used the notation introduced earlier with $\tilde{\psi}_i = \tilde{A}_{m+n\
i}$ and $\tilde{\phi}_i = (-1)^{[i]} \tilde{A}_{i\ m+n}$. We then have
$$
\left( \tilde{P}_k \tilde{P}_{0r} \right)_{ij} = \left( \tilde{P}_k \right)_{i\
m+n}\tilde{\psi}[r]_j \left( \tilde{a}_k - \tilde{a}_{0r} \right)^{-1}.
$$
Multiplication on the left by $\tilde{P}_{0r}$ using equation (\ref{p55star}) gives
\begin{align*}
\left( \tilde{P}_{0r} \tilde{P}_k \tilde{P}_{0r} \right)_{ij} 
&= (-1)^{[i]} \tilde{\phi}[r]_i \left( \tilde{a}_k - q^{2(r)}\tilde{a}_{0r} + (r) q^{(r)} \right)^{-1}
\tilde{\omega}_k \tilde{\psi}[r]_j \left( \tilde{a}_k - \tilde{a}_{0r} \right)^{-1} \\
&= (-1)^{[i]} \tilde{\phi}[r]_i \tilde{\psi}[r]_j \tilde{\omega}_k 
\left( \tilde{a}_k - q^{2(r)} \tilde{a}_{0r} + (r) q^{(r)} \right)^{-1}
\left( \tilde{a}_k - \tilde{a}_{0r} \right)^{-1}.
\end{align*}
Therefore using equation (\ref{p60star}) we obtain
$$
\tilde{P}_{0r} \tilde{P}_k \tilde{P}_{0r} = \tilde{\omega}_{kr} \tilde{P}_{0r},
$$
where now
$$
\tilde{\omega}_{kr} = \tilde{\omega}_k \tilde{\mu}_r\left( \tilde{a}_k - q^{2(r)}
\tilde{a}_{0r} + (r) q^{(r)} \right)^{-1} 
\left( \tilde{a}_k - \tilde{a}_{0r} \right)^{-1}.
$$
\underline{{\em Note}}: Here, if $r\in I_0$ or $k\in I_0$ we obtain zero for the left and
right hand side of the equation. Therefore we only have a non-zero contribution when 
$k\in \overline{I}_0\cup\tilde{I}_1$ and $r\in \overline{I}_0\cup I_1$. In this case
$$
\chi_{(\Lambda,\Lambda_0)}(\tilde{\omega}_{kr}) = 
\left| 
\left\langle 
\left.
\begin{array}{c} \Lambda+\varepsilon_k\\ \Lambda_0+\varepsilon_{0r} \end{array}
\right|
\left|
\begin{array}{c} \varepsilon_1\\ \varepsilon_{1} \end{array}
;
\begin{array}{c} \Lambda\\ \Lambda_0 \end{array}
\right.
\right\rangle
\right|^2.
$$


\subsection{Calculation of reduced Wigner coefficients}\label{Section26} 

Proceeding recursively by the subalgebra chain 
$$
U_q[gl(m|n)] \supset U_q[gl(m|n-1)]\supset \cdots \supset U_q[gl(m|1)]\supset
U_q[gl(m)]\supset\cdots\supset U_q[gl(1)] 
$$
all $U_q[gl(m|n)]$ Wigner coefficients may be expressed as a product of reduced Wigner coefficients
(or isoscalar factors). The non-zero (dual) vector reduced Wigner coefficients have all
been evaluated and are given explicitly by
\begin{align*}
\left\langle 
\left.
\begin{array}{c} \Lambda+\varepsilon_k\\ \Lambda_0 \end{array}
\right|
\left|
\begin{array}{c} \varepsilon_1\\ \dot{0} \end{array}
;
\begin{array}{c} \Lambda\\ \Lambda_0 \end{array}
\right.
\right\rangle &= \tilde{\omega}_k^{\frac12}, 
\\
\left\langle 
\left.
\begin{array}{c} \Lambda-\varepsilon_k\\ \Lambda_0 \end{array}
\right|
\left|
\begin{array}{c} \overline{\varepsilon}_1\\ \dot{0} \end{array}
;
\begin{array}{c} \Lambda\\ \Lambda_0 \end{array}
\right.
\right\rangle &= \omega_k^{\frac12} ,
\\
\left\langle 
\left.
\begin{array}{c} \Lambda+\varepsilon_k\\ \Lambda_0+\varepsilon_{0r} \end{array}
\right|
\left|
\begin{array}{c} \varepsilon_1\\ \varepsilon_1 \end{array}
;
\begin{array}{c} \Lambda\\ \Lambda_0 \end{array}
\right.
\right\rangle &= \tilde{s}_{rk} \tilde{\omega}_k^{\frac12} ,
\\
\left\langle 
\left.
\begin{array}{c} \Lambda-\varepsilon_k\\ \Lambda_0-\varepsilon_{0r} \end{array}
\right|
\left|
\begin{array}{c} \overline{\varepsilon}_1\\ \overline{\varepsilon}_1 \end{array}
;
\begin{array}{c} \Lambda\\ \Lambda_0 \end{array}
\right.
\right\rangle &= s_{rk} \omega_k^{\frac12},
\end{align*}
where the phases $\tilde{s}_{rk}, s_{rk}=\pm 1$ are given by the (known) classical (i.e.
$q\rightarrow 1$ case) phases. Therefore, we have actually calculated, in principle, all
Wigner coefficients.

While the vector representation $\pi$ is unitary, the dual vector representation $\pi^*$
is not. This does not, however, affect the calculation of reduced Wigner coefficients.
Indeed, let $\overline{\pi}$ be the (undeformed) unitary dual representation with
corresponding matrix
$$
\overline{A}_{ij} = q^{-(\rho,\varepsilon_i)} A_{ij} q^{(\rho, \varepsilon_j)}. 
$$ 
Then if $\overline{P}_k$ are the associated projectors we have also
$$
\left( \overline{P}_k \right)_{ij} = q^{(\rho,\varepsilon_j-\varepsilon_i)} 
\left( P_k \right)_{ij}
$$
$$
\Rightarrow \ \ \left( \overline{P}_k \right)_{m+n\ m+n} = \left( P_k \right)_{m+n\ m+n} =
\omega_k,
$$
i.e. the $\omega_k$ indeed give the squared reduced Wigner coefficients for the unitary
case.

Similarly, letting $\overline{\rho}$ denote the $U_q[gl(m|n-1)]$ analogue of the
$U_q[gl(m|n)]$ Weyl vector $\rho$, we have seen that
$$
(\rho-\overline{\rho},\varepsilon_i - \varepsilon_j) = 0, \mbox{ for } i,j<m+n.
$$
Therefore we have, for $i,j<m+n$, in obvious notation, 
\begin{align*}
\left( \overline{P}_{0r}\overline{P}_k\overline{P}_{0r} \right)_{ij}
&= 
\left( \overline{P}_{0r} \right)_{i\ell} \left( \overline{P}_k \right)_{\ell s} \left(
\overline{P}_{0r} \right)_{sj} \\
&=
q^{(\overline{\rho},\varepsilon_\ell - \varepsilon_i)}
q^{(\rho,\varepsilon_s-\varepsilon_\ell)} q^{(\overline{\rho},\varepsilon_j -
\varepsilon_\ell)} 
\left( P_{0r} \right)_{i\ell} \left( P_k \right)_{\ell s} \left( P_{0r} \right)_{sj} \\
&=
q^{(\overline{\rho},\varepsilon_j-\varepsilon_i)} \left( P_{0r} P_k P_{0r} \right)_{ij} \\
&= 
q^{(\overline{\rho},\varepsilon_j-\varepsilon_i)} \omega_{kr} \left( P_{0r} \right)_{ij}
\\
&= \omega_{kr} \left( \overline{P}_{0r} \right)_{ij}.
\end{align*}
Therefore, again, the $\omega_{kr}$ determine the required (unitary) reduced Wigner
coefficients.

Thus the desired reduced Wigner coefficients for the dual vector representation are
independent of the explicit choice of representation $\pi^*$ or $\overline{\pi}$, as might
be expected.


\subsection{Alternative formulae}\label{Section27} 

Given
$$
Y = \frac{1-q^{-2y}}{q-q^{-1}},\ \ X = \frac{1-q^{-2x}}{q-q^{-1}},\ \ (x) = \pm 1,
$$
we note the identity
$$
Y - q^{2(x)}X + (x)q^{(x)} = q^{-x-y+(x)} [y-x+(x)]_q
$$
or
$$
Y - q^{-2(x)}X - (x)q^{-(x)} = q^{-x-y-(x)} [y-x-(x)]_q,
$$
$$
Y-X = q^{-y-x}[y-x]_q.
$$
From this we note the following characteristic root formulae:
\begin{align*}
a_k - q^{-2(r)} a_{0r} - (r) q^{-(r)} 
&= q^{-\alpha_k - \alpha_{0r} - (r)}\left[ \alpha_k - \alpha_{0r} - (r) \right]_q, \\
a_k - a_\ell 
&= q^{-\alpha_k-\alpha_\ell} \left[ \alpha_k - \alpha_\ell \right]_q,\\
\tilde{a}_k - \tilde{a}_\ell
&= q^{-\overline{\alpha}_\ell - \overline{\alpha}_k}\left[
\overline{\alpha}_k-\overline{\alpha}_\ell \right]_q, \\
a_k - a_{0r}
&= q^{-\alpha_k-\alpha_{0r}}\left[ \alpha_k - \alpha_{0r} \right]_q, \\
\tilde{a}_k - \tilde{a}_{0r}
&= q^{ -\overline{\alpha}_{0r} - \overline{\alpha}_k }\left[
\overline{\alpha}_k-\overline{\alpha}_{0r} \right]_q, \\   
a_{0r} - q^{-2(\ell)}a_{0\ell} - (\ell) q^{-(\ell)}
&= q^{ -\alpha_{0r} - \alpha_{0\ell}-(\ell) }\left[ \alpha_{0r} - \alpha_{0\ell} - (\ell) \right]_q, \\
\tilde{a}_{0r} - q^{2(\ell)} \tilde{a}_{0\ell} + (\ell) q^{(\ell)}
&= q^{-\overline{\alpha}_{0r} - \overline{\alpha}_{0\ell}+(\ell) }\left[ \overline{\alpha}_{0r} -
\overline{\alpha}_{0\ell} + (\ell) \right]_q, \\
\tilde{a}_k - q^{2(r)} \tilde{a}_{0r} + (r)q^{(r)}
&= q^{-\overline{\alpha}_k - \overline{\alpha}_{0r} + (r)}\left[ \overline{\alpha}_k -
\overline{\alpha}_{0r} + (r) \right]_q,
\end{align*}
where we recall that
\begin{align*}
a_k &= \frac{1-q^{-2\alpha_k}}{q-q^{-1}},\\
\alpha_i &= \Lambda_i+m-n-i,\\
\alpha_\mu &= \mu - n - \Lambda_\mu,\\
\tilde{a}_k &= \frac{1-q^{-2\overline{\alpha}_k}}{q-q^{-1}}, \\
\overline{\alpha}_i &= \Lambda_i + 1-i,\\
\overline{\alpha}_\mu &= \mu-m-1-\Lambda_\mu,
\end{align*}
with similar expressions for $a_{0r}$ and $\tilde{a}_{0r}.$

Now using
$$
\omega_k = \prod_{r\in I_0\cup I_1} \left(  a_k - q^{-2(r)} a_{0r} - (r) q^{-(r)} \right)
\prod_{\ell\neq k}^{I_0\cup\tilde{I}_1} \left( a_k - a_\ell \right)^{-1}
$$
we obtain
\begin{align*}
\omega_k 
&= 
\frac{ \prod_{r\in I_0\cup I_1} q^{-\alpha_k-\alpha_{0r} - (r)} \left[ \alpha_k -
\alpha_{0r} - (r) \right]_q   }{ \prod_{\ell\neq k}^{I_0\cup \tilde{I}_1} q^{-\alpha_k -
\alpha_\ell} \left[ \alpha_k- \alpha_\ell \right]_q  } \\
&= 
q^{\xi_k} \frac{ \prod_{r\in I_0\cup I_1} \left[ \alpha_k - \alpha_{0r} - (r) \right]_q  }{
\prod_{\ell\neq k}^{I_0\cup \tilde{I}_1} \left[ \alpha_k- \alpha_\ell \right]_q   }
\end{align*}
where the ``q-phase'' $\xi_k$ is given by
\begin{align*}
\xi_k 
&= 
\sum_{\ell\neq k}^{I_0\cup \tilde{I}_0} \alpha_\ell - \sum_{r\in I_0\cup I_1} (\alpha_{0r} + (r)) \\
&=
\sum_{\ell\in I_0\cup \tilde{I}_0 } \alpha_\ell - \sum_{r\in I_0\cup I_1} (\alpha_{0r} + (r)) -
\alpha_k \\
&=
- \left( |I_0| + \alpha_k + \eta(\Lambda, \Lambda_0) \right),
\end{align*}
where we have used the results
$$
\Lambda_i = \Lambda_{0i} + 1, \mbox{ for $i\in I_0$},
$$
$$
\sum_{r\in I_0\cup I_1} (r) = |I_0| - n+1
$$
and where 
$$
\eta(\Lambda,\Lambda_0) \equiv \sum_{\mu=1}^n \Lambda_\mu - \sum_{\mu=1}^{n-1} \Lambda_{0\mu}.
$$
It is interesting to note that $|I_0|+\eta(\Lambda,\Lambda_0)$ gives the eigenvalue of the generator
$E_{m+n\ m+n}$, so that we may write, in a slight abuse of notation,
$$
\xi_k =  -E_{m+n\ m+n} - \alpha_k.
$$
Similarly we have
\begin{align*}
\tilde{\omega}_k 
&=
\prod_{r\in \overline{I}_0\cup I_1} \left( \tilde{a}_k - q^{2(r)} \tilde{a}_{0r} + (r)q^{(r)}\right) 
\prod_{\ell\neq k}^{\overline{I}_0\cup \tilde{I}_1} \left( \tilde{a}_k - \tilde{a}_\ell \right)^{-1}
\\
&=
\frac{ \prod_{r}^{\overline{I}_0\cup I_1} q^{-\overline{\alpha}_k - \overline{\alpha}_{0r} + (r)}
\left[ \overline{\alpha}_k - \overline{\alpha}_{0r} + (r) \right]_q  }{ \prod_{\ell\neq
k}^{\overline{I}_0\cup \tilde{I}_1} q^{-\overline{\alpha}_k - \overline{\alpha}_\ell} 
\left[ \overline{\alpha}_k - \overline{\alpha}_\ell \right]_q  }
\\
&=
q^{\tilde{\xi}_k} \frac{ \prod_r^{\overline{I}_0\cup I_1} \left[ \overline{\alpha}_k -
\overline{\alpha}_{0r} + (r) \right]_q  }{ \prod_{\ell\neq k}^{\overline{I}_0\cup\tilde{I}_1} 
\left[ \overline{\alpha}_k - \overline{\alpha}_\ell \right]_q},
\end{align*}
where
\begin{align*}
\tilde{\xi}_k 
&=
\sum_{\ell\neq k}^{\overline{I}_0\cup \tilde{I}_1} \overline{\alpha}_\ell -
\sum_{r}^{\overline{I}_0\cup I_1} \left( \overline{\alpha}_{0r} - (r) \right)
\\
&= 
\sum_\ell^{\overline{I}_0\cup \tilde{I}_1} \overline{\alpha}_\ell - \sum_r^{\overline{I}_0\cup I_1}
\overline{\alpha}_{0r} + \sum_r^{\overline{I}_0\cup I_1} - \overline{\alpha}_k.
\end{align*}
Using
$$
\sum_r^{\overline{I}_0\cup I_1} (r) = |\overline{I}_0 | - |I_1| = |\overline{I}_0| - n+1 = m-n+1 -
|I_0|,
$$
$$
\sum_r^{\overline{I}_0} \left( \overline{\alpha}_r - \overline{\alpha}_{0r} \right)
= \sum_{r}^{\overline{I}_0} \left( \Lambda_r - \Lambda_{0r} \right) = 0,
$$
$$
\sum_{\mu=1}^n \overline{\alpha}_\mu - \sum_{\mu=1}^{n-1}\overline{\alpha}_{0\mu} = 
n-m-1-\eta(\Lambda,\Lambda_0),
$$
gives
\begin{align*}
\tilde{\xi}_k &=
-\left( \eta(\Lambda,\Lambda_0) + |I_0| + \overline{\alpha}_k \right) \\
&= (\nu,\delta_{m+n}) - \overline{\alpha}_k,
\end{align*}
with $\nu$ being the weight of the Gelfand-Tsetlin state on which the operator acts.

For the squared reduced matrix elements, we have, from (\ref{p58star}),
\begin{align*}
\mu_r 
&=
(-1)^{|I_0| + n-1} \prod_k^{I_0 \cup \tilde{I}_1} (a_k-a_{0r})\prod_{\ell\neq r}^{I_0\cup I_1} 
\left( a_{0r} - q^{-2(\ell)}a_{0\ell} - (\ell)q^{-(\ell)} \right)^{-1} \\
&= 
(-1)^{|I_0| + n-1} \frac{ \prod_k^{I_0 \cup \tilde{I}_1} q^{-\alpha_k - \alpha_{0r}}
[\alpha_k-\alpha_{0r}]_q }{ \prod_{\ell\neq r}^{I_0\cup I_1} q^{-\alpha_{0r}-\alpha_{0\ell} -
(\ell)} [\alpha_{0r}-\alpha_{0\ell} - (\ell)]_q } \\
&= 
(-1)^{|I_0| + n-1} q^{\eta_r} \frac{ \prod_k^{I_0 \cup \tilde{I}_1} [\alpha_k-\alpha_{0r}]_q }{
\prod_{\ell\neq r}^{I_0\cup I_1} [\alpha_{0r}-\alpha_{0\ell} - (\ell)]_q  },
\end{align*}
where
\begin{align*}
\eta_r 
&= 
\sum_{\ell\neq r}^{I_0\cup I_1} (\alpha_{0\ell} + (\ell)) + \sum_{\ell\neq r}^{I_0\cup
I_1}\alpha_{0r} - \sum_k^{I_0\cup\tilde{I}_1} \alpha_{0r} - \sum_k^{I_0\cup \tilde{I}_1} \alpha_k \\
&= 
\sum_\ell^{I_0\cup I_1}(\alpha_{0\ell} + (\ell)) - (\alpha_{0r}+(r)) - 2\alpha_{0r} -
\sum_k^{I_0\cup \tilde{I}_1 \alpha_k} \\
&= \sum_\ell^{I_0\cup I_1} \alpha_{0\ell} - \sum_k^{I_0\cup \tilde{I}_1} \alpha_k +
\sum_\ell^{I_0\cup I_1} (\ell) - 3\alpha_{0r} -(r).
\end{align*}
Using 
\begin{align*}
\sum_\ell^{I_0\cup I_1} (\ell) &= |I_0| - n+1,\\
\sum_{\ell}^{I_0\cup I_1} \alpha_{0\ell} - \sum_k^{I_0\cup \tilde{I}_1} \alpha_k
&=
\sum_r^{I_0}( \alpha_{0r} - \alpha_r ) + \sum_\mu\alpha_{0\mu} - \sum_\mu \alpha_\mu\\
&=
\eta(\Lambda,\Lambda_0)+n-1,
\end{align*}
we obtain
$$
\eta_r = |I_0| + \eta(\Lambda,\Lambda_0) - 3\alpha_{0r}-(r) = E_{m+n\ m+n} - 3\alpha_{0r}-(r).
$$
Finally for $\tilde{\mu}_r$ we have
\begin{align*}
\tilde{\mu}_r 
&=
(-1)^{|\overline{I}_0| + n-1} \prod_k^{\overline{I}_0\cup\tilde{I}_1} (\tilde{a}_k-\tilde{a}_{0r})
\prod_{\ell\neq r}^{\overline{I}_0\cup I_1} 
\left( \tilde{a}_{0r} - q^{2(\ell)} \tilde{a}_{0\ell} + (\ell)q^{(\ell)} \right)^{-1} \\
&=
(-1)^{|\overline{I}_0| + n-1}
\frac{ \prod_k^{\overline{I}_0\cup\tilde{I}_1} q^{-\overline{\alpha}_k-\overline{\alpha}_{0r}}
[\overline{\alpha}_k-\overline{\alpha}_{0r}]_q }{ \prod_{\ell\neq r}^{\overline{I}_0\cup I_1}
q^{-\overline{\alpha}_{0r}-\overline{\alpha}_{0\ell} + (\ell) } [
\overline{\alpha}_{0r}-\overline{\alpha}_{0\ell} + (\ell) ]_q } \\
&=
(-1)^{|\overline{I}_0| + n-1}
q^{\tilde{\eta}_r} \frac{ \prod_k^{\overline{I}_0\cup\tilde{I}_1}
[\overline{\alpha}_k-\overline{\alpha}_{0r}]_q  }{ \prod_{\ell\neq r}^{\overline{I}_0\cup I_1}[
\overline{\alpha}_{0r}-\overline{\alpha}_{0\ell} + (\ell) ]_q  },
\end{align*}
where
\begin{align*}
\tilde{\eta}_r 
&=
\sum_{\ell\neq r}^{\overline{I}_0\cup I_1} ( \overline{\alpha}_{0\ell} - (\ell) ) -
\sum_k^{\overline{I}_0\cup\tilde{I}_1} \overline{\alpha}_k +\sum_{\ell\neq r}^{\overline{I}_0\cup
I_1} \overline{\alpha}_{0r} - \sum_k^{\overline{I}_0\cup\tilde{I}_1} \overline{\alpha}_{0r} \\
&=
\sum_{\ell}^{\overline{I}_0\cup I_1} \overline{\alpha}_{0\ell} - \sum_k^{\overline{I}_0\cup\tilde{I}_1} \overline{\alpha}_k 
-\overline{\alpha}_{0r} + (r) - 2\overline{\alpha}_{0r} - \sum_\ell^{\overline{I}_0\cup I_1} (\ell).
\end{align*}
Now
\begin{align*}
\sum_{\ell}^{\overline{I}_0\cup I_1}  (\ell)
&= 
|\overline{I}_0| - (n-1) = m-n+1-|I_0|,\\
\sum_\ell \overline{\alpha}_{0\ell} - \sum_k \overline{\alpha}_k
&=
\sum_{\ell\in \overline{I}_0} (\overline{\alpha}_{0\ell} - \overline{\alpha}_\ell) +
\sum_{\mu=1}^{n-1} \overline{\alpha}_{0\mu}- \sum_{\mu=1}^n \overline{\alpha}_\mu \\
&=
\eta(\Lambda,\Lambda_0) + 1+m-n.
\end{align*}
We therefore obtain
\begin{align*}
\tilde{\eta}_r &= \eta(\Lambda,\Lambda_0) + |I_0| - 3\overline{\alpha}_{0r} + (r)\\
&= -3\overline{\alpha}_{0r} + (r) - (\nu,\delta_{m+n}),
\end{align*}
with $\nu$ being the weight of the Gelfand-Tsetlin state, as before.


\subsection{Alternative formulae for $\omega_{kr}$, $\tilde{\omega}_{kr}$}\label{Section28} 

Using results above we have
\begin{align*}
\omega_{kr} 
=&
\omega_k\mu_r (a_k - q^{-2(r)}a_{0r} - (r) q^{-(r)})^{-1} (a_k - a_{0r})^{-1} \\
=&
(-1)^{|I_0|+n-1} q^{\xi_k + \eta_r}
\, \frac{ \prod_{s}^{I_0\cup I_1}[\alpha_k-\alpha_{0s}-(s)]_q }{ \prod_{\ell\neq k}^{I_0\cup
\tilde{I}_1} [\alpha_k-\alpha_{\ell}]_q }
\, \frac{ \prod_p^{I_0\cup\tilde{I}_1} [\alpha_p-\alpha_{0r}]_q }{ \prod_{j\neq r}^{I_0\cup I_1}
[\alpha_{0r}-\alpha_{0j}-(j)]_q }
\\
&\quad
q^{\alpha_k+\alpha_{0r}+(r)} 
[\alpha_k-\alpha_{0r}-(r)]_q^{-1}
q^{\alpha_k+\alpha_{0r}}
[\alpha_k-\alpha_{0r}]_q^{-1}  \\
=& 
(-1)^{|I_0|+n-1} q^{\xi_{kr}} 
\prod_{\ell\neq r}^{I_0\cup I_1} \frac{ [\alpha_k-\alpha_{0\ell} - (\ell)]_q }{ [\alpha_{0r} -
\alpha_{0\ell} - (\ell)]_q }
\prod_{p\neq k}^{I_0\cup \tilde{I}_1} \frac{ [\alpha_p - \alpha_{0r}]_q }{ [\alpha_k-\alpha_p]_q }
\\
=&
q^{\xi_{kr}} \prod_{\ell\neq r}^{I_0\cup I_1} \frac{ [\alpha_k-\alpha_{0\ell}-(\ell)]_q }{
[\alpha_{0r} - \alpha_{0\ell} - (\ell)]_q } \prod_{p\neq k}^{I_0\cup \tilde{I}_1} \frac{
[\alpha_p-\alpha_{0r}]_q }{ [\alpha_p-\alpha_k]_q },
\end{align*}
where
\begin{align*}
\xi_{kr} 
&= 
\xi_k + \eta_r + 2(\alpha_k+\alpha_{0r})+(r) \\
&=
|I_0| + \eta(\Lambda,\Lambda_0) - 3\alpha_{0r} - (r) - |I_0| - \alpha_k - \eta(\Lambda,\Lambda_0) +
2(\alpha_k + \alpha_{0r}) + (r)\\
&=
\alpha_k - \alpha_{0r},
\end{align*}
(c.f. the $U_q[gl(n)]$ case).
Similarly
\begin{align*}
\tilde{\omega}_{kr} 
=& 
\tilde{\omega}_k \tilde{\mu}_r \left( \tilde{a}_k - q^{2(r)} \tilde{a}_{0r} + (r) q^{(r)} \right)^{-1}
\left( \tilde{a}_k - \tilde{a}_{0r} \right)^{-1} \\
=& 
q^{\tilde{\xi}_{kr}} \prod_{\ell\neq r}^{\overline{I}_0\cup I_1} \frac{
[\overline{\alpha}_k-\overline{\alpha}_{0\ell}+(\ell)]_q }{
[\overline{\alpha}_{0r} - \overline{\alpha}_{0\ell} + (\ell)]_q } \prod_{p\neq k}^{\overline{I}_0\cup \tilde{I}_1} \frac{
[\overline{\alpha}_p-\overline{\alpha}_{0r}]_q }{ [\overline{\alpha}_p-\overline{\alpha}_k]_q },
\end{align*}
where
\begin{align*}
\tilde{\xi}_{kr} &= \tilde{\xi}_k + \tilde{\eta}_r + 2(\overline{\alpha}_k + \overline{\alpha}_{0r})
- (r)\\
&= -\eta(\Lambda,\Lambda_0) - |I_0| - \overline{\alpha}_k + \eta(\Lambda,\Lambda_0) + |I_0| -
3\overline{\alpha}_{0r} + (r) + 2(\overline{\alpha}_k + \overline{\alpha}_{0r}) - (r)\\
&= \overline{\alpha}_k - \overline{\alpha}_{0r}.
\end{align*}
It is worth pointing out the remarkable similarity between the formulae derived above, and those of
$U_q[gl(n)]$.



\section{Concluding remarks}

We have established formulae for certain invariants of $U_q[gl(m|n)]$, including the reduced Wigner
coefficients and reduced matrix elements. These results will be utilised in a forthcoming paper on
the matrix elements of the quantum supergroup generators.



\appendix

\section{Appendix: The $u$-operator}\label{AppendixA}

Here we present the proofs of the results relating to the $u$-operator defined in the main body of the
article in equation (\ref{uoperator}).

Following the notation in the text, with $R=a_i\otimes b_i$ (sum on $i$), we write
$$
u = (-1)^{[a_i]}S(b_i)a_i.
$$
We have, from equation (\ref{p1starstar}),
\begin{align*}
	& \ \ (-1)^{[b_i][a_{(1)}]}a_ia_{(1)}\otimes  b_ia_{(2)} =
	(-1)^{[a_{(1)}]([a_{(2)}]+[a_i])}a_{(2)}a_i\otimes a_{(1)}b_i\\
	\stackrel{T}{\Rightarrow} & \ \ (-1)^{([b_i]+[a_{(2)}])([a_i]+[a_{(1)}]) +
	[b_i][a_{(1)}]}b_ia_{(2)}\otimes a_ia_{(1)} \\
	& \ \ \ \ =
	(-1)^{[a_{(1)}]([a_{(2)}]+[a_i])+([a_{(2)}]+[a_i])([a_{(1)}]+[b_i])}
	a_{(1)}b_i\otimes a_{(2)}a_i\\
	\Rightarrow  & \ \  (-1)^{[b_i] +
	[a_{(2)}]([a_i]+[a_{(1)}])}b_ia_{(2)}\otimes a_ia_{(1)} =
	(-1)^{[b_i]([a_{(2)}]+[a_i])}a_{(1)}b_i\otimes a_{(2)}a_i\\
	\stackrel{S\otimes \id}{\Rightarrow} & \ \
	(-1)^{[b_i]+[a_{(1)}][a_{(2)}]}S(a_{(2)})S(b_i)\otimes a_ia_{(1)} =
	(-1)^{[b_i]([a]+[a_i])} S(b_i)S(a_{(1)})\otimes a_{(2)}a_i
\end{align*}
\begin{equation}
	\stackrel{m}{\Rightarrow} \ \
	(-1)^{[a_{(1)}][a_{(2)}]}S(a_{(2)})ua_{(1)} = \varepsilon(a)u
	\label{p3star}
\end{equation}
Now set 
$$
R^{-1} = c_i\otimes d_i.
$$
We prove the following result.
\begin{lemma} (Notation as above):\\ 
(i) $S^2(a)u = ua,$ $\forall a\in H$ \\
(ii) $u$ is invertible with $u^{-1} =
			(-1)^{[c_i]}S^{-1}(d_i)c_i.$
\end{lemma}
Proof: (i) First observe 
\begin{align*}
	(S\otimes S)\Delta^T(a) &= \Delta(S(a)) \\
	\Rightarrow \ \ (-1)^{[a_{(1)}][a_{(2)}]}S(a_{(2)})\otimes S(a_{(1)}) &=
	S(a)_{(1)}\otimes S(a)_{(2)}.
\end{align*}
Therefore, from (\ref{p3star}), with the element $a$ replaced by $S(a)$, we
have
$$
S^2(a_{(1)})uS(a_{(2)}) = \varepsilon(S(a))u = \varepsilon(a)u
$$
\begin{align*}
\Rightarrow \ \ ua &= \varepsilon(a_{(1)})ua_{(2)}\\
&= S^2(a_{(1)})uS(a_{(2)})a_{(3)} = S^2(a)u.
\end{align*}
(ii)
\begin{align*}
	u\cdot u^{-1} &= u\cdot (-1)^{[c_j]}S^{-1}(d_j)c_j\\
	& \stackrel{(i)}{=} (-1)^{[c_j]} S(d_j)uc_j\\
&= (-1)^{[c_j]+[a_i]}S(d_j)S(b_i)a_ic_j\\
&=  (-1)^{[a_i]+[c_j]+[a_i][c_j]}S(b_id_j)a_ic_j\\
&= m(S\otimes\id)[(b_i\otimes a_i)(d_j\otimes c_j)](-1)^{[a_i]+[c_j]}\\
&= m(S\otimes\id)[R^T(R^{-1})^T] = I.
\end{align*}

\begin{flushright}$\Box$\end{flushright}

\begin{cor}
$S^2(a) = uau^{-1},\ \ \forall a\in H.$
\end{cor}

Our aim here is to prove the following.

\begin{thm} \label{Theorem1}
$\Delta(u) = (u\otimes u)(R^TR)^{-1}.$
\end{thm}

Before proceeding, note that
$$
(\Delta\otimes\id)R = R_{13}R_{23},\ \ (\id\otimes R) = R_{13}R_{12} 
$$
\begin{align*}
\Rightarrow \ \ (\Delta^T\otimes \id)R &= (T\otimes\id)(R_{13}R_{23})\\
&= (T\otimes \id)(a_i\otimes a_j\otimes b_ib_j)(-1)^{[b_i][a_j]}\\
&= (a_j\otimes a_i\otimes b_ib_j)(-1)^{[b_i][a_j] + [a_i][a_j]}\\
&= (I\otimes a_i\otimes b_i)(a_j\otimes I\otimes b_j)\\
&= R_{23}R_{13}
\end{align*}
and similarly
$$
(\id\otimes \Delta^T)R = R_{12}R_{13}.
$$
An immediate consequence of this -- one of the key results of Drinfeld \cite{D86} -- is
that
\begin{align}
R_{12}R_{13}R_{23} &= R_{12}[(\Delta\otimes\id)R] \nonumber\\
&= [(\Delta^T\otimes \id)R] R_{12}\nonumber\\
&= R_{23}R_{13}R_{12}, \label{qybe}
\end{align}
known as the quantum Yang-Baxter equation.
Now observe that 
\begin{align*}
\Delta(u) R^TR &= (-1)^{[a_i]} \Delta(S(b_i))\Delta(a_i)R^TR\\
&= (-1)^{[b_i]}(S\otimes S)\Delta^T(b_i)R^TR\Delta(a_i)\\
&= \varphi\left( R_{12}^TR_{12}(\Delta^T\otimes \Delta^T)R \right),
\end{align*}
where $\varphi:H^{\otimes 4}\longrightarrow H^{\otimes 2}$ is defined by
\begin{align*}
\varphi(c_1\otimes c_2\otimes c_3\otimes c_4)
&= (-1)^{ ( [c_1]+[c_2] ) [c_3] + [c_2][c_4] } S(c_3)c_1\otimes S(c_4)c_2\\
&= (-1)^{\gamma} (S(c_3)\otimes S(c_4))(c_1\otimes c_2),
\end{align*}
where $\gamma = ([c_1]+[c_2])([c_3]+[c_4])$, since
\begin{align*}
\varphi\left( R_{12}^TR_{12}(\Delta\otimes\Delta^T)R \right) 
&= \varphi\left(  R_{12}^TR_{12}(\Delta(a_i)\otimes\Delta^T(b_i)) \right)\\
&= \varphi\left( R_{12}^TR_{12}(\Delta(a_i)\otimes b_{i(1)}\otimes b_{i(2)}) \right)
(-1)^{ [b_{i(1)}][b_{i(2)}] }\\
&= (S(b_{i(2)})\otimes S(b_{i(1)}))R^TR\Delta(a_i)(-1)^{ [b_{i(1)}][b_{i(2)}] + [a_i][b_i]
} \\
&= (S\otimes S)\Delta^T(b_i)R^TR\Delta(a_i)(-1)^{[a_i]}\\
&= \Delta(u)R^TR.
\end{align*}
Now observe that
\begin{align*}
(\Delta\otimes\Delta^T)R &= (\Delta\otimes\id\otimes\id)(\id\otimes \Delta^T)R\\
&= (\Delta\otimes \id\otimes\id)R_{12}R_{13}\\
&= R_{13}R_{23}R_{14}R_{24}.
\end{align*}
Therefore,
\begin{align}
\Delta(u)R^TR &= \varphi\left( R_{12}^TR_{12}R_{13}R_{23}R_{14}R_{24} \right) \nonumber \\
& \stackrel{(\ref{qybe})}{=} \varphi \left( R_{12}^TR_{23}R_{13}R_{12}R_{14}R_{24}
\right).
\label{p5star}
\end{align} 
Consider the following Lemma.

\begin{lemma}
For all $c\in H^{\otimes 4}$, we have
\begin{itemize}
\item[(i)]
$\varphi\left\{ R_{12}^TR_{23}c \right\} = \varphi(c),$
\item[(ii)]
$\varphi\left\{ R_{13}c \right\} = (u\otimes I)\overline{\varphi}(c),$
where 
$$
\overline{\varphi}(c) = (-1)^\gamma (S^{-1}(c_3)\otimes S(c_4))(c_1\otimes c_2),\ \
\gamma=([c_1]+[c_2])([c_3]+[c_4]).
$$
\end{itemize}
\end{lemma}
Proof:
(i) First note, in obvious notation, that
\begin{align*}
R_{12}^TR_{23}c &= (b_i\otimes a_i\otimes I\otimes I)(I\otimes a_j\otimes b_j\otimes I)(c_1\otimes
c_2\otimes c_3\otimes c_4)(-1)^{[a_i]}\\
&=(-1)^{\gamma+[a_i]}(b_i\otimes a_i\otimes I\otimes I)(I\otimes a_j\otimes b_j\otimes I)
(I\otimes I\otimes c_3\otimes c_4)(c_1\otimes c_2\otimes I\otimes I),
\end{align*}
where $\gamma=([c_1]+[c_2])([c_3]+[c_4])$,
$$
\Rightarrow \ \ R_{12}^TR_{23}c = (-1)^{\gamma+[a_i]}(b_i\otimes a_ia_j\otimes b_jc_3\otimes
c_4)(c_1\otimes c_2\otimes I\otimes I).
$$
Therefore 
\begin{align*}
\varphi(R_{12}^TR_{23}c) &= (-1)^{\gamma+[a_i]+([c_3]+[c_4]+[b_j])[a_j]}( S(b_jc_3)\otimes S(c_4) )(
b_i\otimes a_i a_j )( c_1\otimes c_2 ) \\
&= (-1)^{\gamma+[a_i]+([c_3]+[c_4]+[b_j])[a_j]}( S(c_3)S(b_j)\otimes S(c_4) )( b_i\otimes a_ia_j )(
c_1\otimes c_2 ) \\
&= (-1)^{\gamma+[a_i]+[a_j]}( S(c_3)\otimes S(c_4) )( S(b_j)b_i\otimes a_ia_j )( c_1\otimes c_2 ).
\end{align*}
Now observe
\begin{align*}
(-1)^{[a_i]+[a_j]} S(b_j)b_i\otimes a_ia_j &= (-1)^{[a_i]+[a_j]+[a_i][a_j]}
(S\otimes\id)(S^{-1}(b_i)b_j\otimes a_ia_j)\\
&= (-1)^{[a_i]+[a_j]} (S\otimes\id) \left[ (S^{-1}(b_i)\otimes a_i)(b_j\otimes a_j) \right] \\
&= (S\otimes\id) \left[ (S^{-1}\otimes\id)R^T\cdot R^T\right]\\
&= (S\otimes\id) \left[ \left(R^T\right)^{-1}\cdot R^T\right]\\
&= I\otimes I.
\end{align*}
$$
\Rightarrow \ \ \varphi(R_{12}^TR_{23}c) = (-1)^\gamma(S(c_3)\otimes S(c_4))(c_1\otimes c_2) =
\varphi(c).
$$
(ii) In the above notation we have
\begin{align*}
R_{13}c &= (a_i\otimes I \otimes b_i\otimes I)(c_1\otimes c_2\otimes c_3\otimes c_4)\\
&= (-1)^{[b_i]([c_1]+[c_2])} a_ic_1\otimes c_2\otimes b_ic_3\otimes c_4\\
\Rightarrow \ \ \varphi(R_{13}c) &=
(-1)^{[b_i]([c_1]+[c_2])+([b_i]+[c_3]+[c_4])([a_i]+[c_1]+[c_2])} (S(b_ic_3)\otimes
S(c_4))(a_ic_1\otimes c_2) \\
&= (-1)^{\gamma+[a_i]([b_i]+[c_3]+[c_4])} (S(b_ic_3)\otimes S(c_4))(a_ic_1\otimes c_2), 
\end{align*}
where $\gamma=([c_1]+[c_2])([c_3)+[c_4])$ as above. Thus
\begin{align*}
\varphi(R_{13}c) &= (-1)^{\gamma+[a_i]([b_i]+[c_4])} (S(c_3)S(b_i)\otimes
S(c_4))(a_ic_1\otimes c_2)\\
&= (-1)^{\gamma+[a_i]} (S(c_3)\otimes S(c_4))(S(b_i)a_ic_1\otimes c_2)\\
&= (-1)^{\gamma}(S(c_3)\otimes S(c_4))(uc_1\otimes c_2)\\
&= (-1)^{\gamma} (S(c_3)u\otimes S(c_4))(c_1\otimes c_2)\\
&= (-1)^{\gamma}(u\otimes I)(S^{-1}(c_3)\otimes S(c_4))(c_1\otimes c_2)\\
&= (u\otimes I)\overline{\varphi}(c).
\end{align*}

\begin{flushright}$\Box$\end{flushright}

By repeated application of this lemma, we obtain
\begin{align*}
\Delta(u)R^TR &\stackrel{(\ref{p5star})}{=} \varphi\left(
R_{12}^TR_{23}R_{13}R_{12}R_{14}R_{24} \right) \\
&\stackrel{(i)}{=} \varphi\left( R_{13}R_{12}R_{14}R_{24} \right)\\ 
&\stackrel{(ii)}{=} (u\otimes I) \overline{\varphi}\left( R_{12}R_{14}R_{24} \right)\\
&= (u\otimes I) \overline{\varphi}\left( a_ia_j\otimes b_ia_k\otimes I\otimes b_jb_k \right)
(-1)^{[a_j]([a_i]+[a_k])} \\
&= (u\otimes I)(I\otimes S(b_jb_k))(a_ia_j\otimes b_ia_k)(-1)^{[a_j]([a_i]+[a_k]) +
([b_j]+[b_k])([a_j]+[a_k])} \\
&= (u\otimes I)(-1)^{[a_j]([a_i]+[a_k])+[a_j]+[a_k]}(I\otimes S(b_jb_k))(a_ia_j\otimes
b_ia_k)\\
&= (u\otimes I) (-1)^{[a_i][a_j]+[a_j]+[a_k]} (I\otimes S(b_k)S(b_j))(a_ia_j\otimes
b_ia_k)\\
&= (u\otimes I)(I\otimes S(b_k))(I\otimes S(b_j))(a_ia_j\otimes b_i)(I\otimes
a_k)(-1)^{[a_i][a_j]+[a_j]+[a_k]} \\
&= (u\otimes I)(I\otimes S(b_k))(I\otimes S(b_j))(a_ia_j\otimes S(b_j)b_i)(I\otimes
a_k)(-1)^{[a_k]}
\end{align*}
where
\begin{align*}
a_ia_j\otimes S(b_j)b_i &= (\id\otimes S)(a_ia_j\otimes S^{-1}(b_i)b_j)(-1)^{[a_i][a_j]}\\
&= (\id\otimes S)\left( (a_i\otimes S^{-1}(b_i))(a_j\otimes b_j) \right) \\
&= (\id\otimes S)\left( R^{-1}R \right) = I\otimes I \\
\Rightarrow \ \ \Delta(u)R^TR &= (u\otimes I) (I\otimes S(b_k)a_k)(-1)^{[a_k]}\\
&= (u\otimes I)(I\otimes u) = (u\otimes u).
\end{align*}
This in fact completes the proof of Theorem \ref{Theorem1}.

\noindent
\underline{{\em Remark}:} Since $R^{-1}=c_i\otimes d_i = (S\otimes\id)R=S(a_i)\otimes b_i$
we must have
\begin{align}
u^{-1} &= S^{-1}(d_i)c_i(-1)^{[c_i]} \nonumber\\
&= S^{-1}(b_i)S(a_i)(-1)^{[a_i]} \nonumber\\
&= (-1)^{[a_i]}S^{-2}(b_i)a_i \label{p8star}
\end{align} 
which follows from $R=(S\otimes S)R = (S^{-1}\otimes S^{-1})R.$

We now observe that 
$$
\tilde{R} = (R^T)^{-1}
$$
is also an $R$-matrix. We denote the corresponding
$u$-operator by $\tilde{u}$. We have (c.f. the non-super case)

\begin{lemma}
$\tilde{u} = S(u^{-1}).$
\end{lemma}
Proof: 
\begin{align*}
\tilde{R} &= \tilde{a}_i\otimes\tilde{b}_i\\
&= (R^{-1})^T\\
&= ((S\otimes\id)R)^T \\
&= (S(a_i)\otimes b_i)^T \\
&= (-1)^{[a_i]} b_i\otimes S(a_i).
\end{align*}
Therefore, for the $u$-operator of $\tilde{R}$ we have
\begin{align*}
\tilde{u} &= (-1)^{[\tilde{a}_i]}S(\tilde{b}_i)\tilde{a}_i\\
&= S^2(a_i)b_i\\
&= S\left( S^{-1}(b_i)S(a_i) \right) (-1)^{[a_i]}\\
&= S(S^{-2}(b_i)a_i(-1)^{[a_i]})\\
&\stackrel{(\ref{p8star})}{=} S(u^{-1}).
\end{align*}

\begin{flushright}$\Box$\end{flushright}

Hence we have
\begin{align*}
\Delta(\tilde{u}) &= (\tilde{u}\otimes\tilde{u}) (\tilde{R}^T\tilde{R})^{-1}\\
&= (\tilde{u}\otimes\tilde{u})(R^{-1}(R^T)^{-1})^{-1}\\
&= (\tilde{u}\otimes \tilde{u})(R^TR).
\end{align*}

\begin{cor} \label{Corollary2}
$u\otimes u$ commutes with $R^TR$.
\end{cor}
Proof: Follows since $\Delta(u)$ commutes with $R^TR$.

\begin{cor}
$\Delta(S(u)) = (S(u)\otimes S(u))(R^TR)^{-1}$.
\end{cor}
Proof: Follows from Corollary \ref{Corollary2} and the fact that $S(u)$ is the inverse of
$\tilde{u} = S(u^{-1}).$


\section{Appendix: The $L$-operator}\label{AppendixB}

For completeness, here we summarise some noteworthy properties of the $L$-operator under the action
of the co-product. This provides a powerful check on the $L$-operator.

From equation (\ref{p25starstar}) we have
$$
(\pi\otimes\id)R = \sum_{i\leq j} e_{ji}\otimes \tilde{E}_{ij},
$$
where 
$$
\tilde{E}_{ij} = \left\{ 
\begin{array}{rl}
(q-q^{-1})(-1)^{[i]}q^{\frac12((i)E_{ii} + (j)E_{jj} - (i))}E_{ij} &,\ i<j,\\
q^{(i)E_{ii}} &,\ i=j,
\end{array}
\right.
$$
where $E_{ij}$ is defined recursively by
$$
E_{ij} = E_{ik}E_{kj} - q^{-(k)}E_{kj}E_{ik}, \ \ i\lessgtr k \lessgtr j, 
$$
with 
$$
E_{i\ i+1} = e_i,\ \ E_{i+1\ i} = f_i.
$$
It is suggestive that we write
$$
(\pi\otimes\id)R = \sum_i e_{ii}\otimes q^{(i)E_{ii}} + (q-q^{-1})\sum_{i< j} e_{ji}\otimes
\overline{E}_{ij},
$$
where
$$
\overline{E}_{ij} = (-1)^{[i]} q^{\frac12 ((i)E_{ii} + (j)E_{jj} - (i))} E_{ij},\ \ i\neq j.
$$
Then using induction with 
$$
\Delta(e_i) = q^{h_i/2}\otimes e_i + e_i\otimes q^{-h_i/2}
$$
we arrive at the coproduct rule
\begin{align*}
\Delta(E_{ij}) =& q^{h_{ij}/2}\otimes E_{ij} + E_{ij}\otimes q^{-h_{ij}/2} \\ 
& \quad +(q-q^{-1})\sum_{i<k<j} (-1)^{[k]} q^{-(k)} 
\left( q^{-(j)E_{jj}/2} E_{ik} q^{(k)E_{kk}/2} \right)
\otimes \left( q^{(k)E_{kk}/2} E_{kj} q^{-(i)E_{ii}/2} \right),\ \ i<j,
\end{align*}
where
$$
h_{ij} = (i)E_{ii} - (j)E_{jj}.
$$

From this result we deduce the following action for the coproduct on $\overline{E}_{ij}$:
\begin{align*}
\Delta(\overline{E}_{ij}) 
=& 
(-1)^{[i]} q^{-(i)/2} 
\Delta( q^{(i)E_{ii}/2 + (j)E_{jj}/2} E_{ij} ) \\
=& 
(-1)^{[i]} q^{-(i)/2} 
\left\{ \left( q^{(i)E_{ii}/2 + (j)E_{jj}/2} \right)\otimes \left( q^{(i)E_{ii}/2 + (j)E_{jj}/2}
\right) \right\} \Delta(E_{ij}) \\
=&
q^{(i)E_{ii}}\otimes \overline{E}_{ij} + \overline{E}_{ij}\otimes q^{(j)E_{jj}}
+ (q-q^{-1}) \sum_{i<k<j}\overline{E}_{ik}\otimes \overline{E}_{kj},\ \ i<j. 
\end{align*}
This implies the following remarkably simple coproduct formula for the $\tilde{E}_{ij}$:
\begin{align*}
\Delta(\tilde{E}_{ij}) &= q^{(i)E_{ii}} \otimes \tilde{E}_{ij} + \tilde{E}_{ij}\otimes q^{(j)E_{jj}}
+ \sum_{i<k<j} \tilde{E}_{ik}\otimes \tilde{E}_{kj},\ \ i<j,\\
\Delta(\tilde{E}_{ii}) &= \Delta(q^{(i)E_{ii}})  = \tilde{E}_{ii}\otimes \tilde{E}_{ii}.
\end{align*}
Using the above coproduct rules it is easy to verify the identity
$$
(\pi\otimes\id\otimes\id)(\id\otimes\Delta)R = (\pi\otimes\id\otimes \id) R_{13}R_{12}
$$
as required of the $R$-matrix. This then provides a powerful check on the formula for the
$L$-operator.

%
%
%

%

\end{document}